\documentclass[jip]{degruyter-journal-a}      
 \firstpage{1}
 \volume{}
 \pubyear{2014}
 \aopyear{}
 \doiyear{}
 \doi{}
 \communicated{}
 \received{November 23, 2013}
 \revised{February 18, 2014}
 \accepted{February 20, 2014}
 
 \accepted{}

\usepackage{graphicx}
\usepackage{graphicx}
\usepackage{amssymb, amsmath}
\usepackage{bm}
\usepackage{mathrsfs}

\newtheorem{defi}{Definition}[section]
\newtheorem{lem}{Lemma}[section]
\newtheorem{vet}{Theorem}[section]
\newtheorem{dos}{Corollary}[section]
\newtheorem{tvr}{Proposition}[section]
\newtheorem{ex}{Example}[section]
\newtheorem{rem}{Remark}[section]
\newcommand{\R}{\mathbb{R}}
\newcommand{\CC}{\mathbb{C}}
\newcommand{\Z}{\mathbb{Z}}
\newcommand{\N}{\mathbb{N}}
\newcommand{\K}{\mathcal{K}}

\newcommand{\vecx}{\bm{x}}

\newcommand{\vecn}{\bm{n}}
\newcommand{\vect}{\bm{t}}
\newcommand{\ds}{\hbox{d}s}
\newcommand{\du}{\hbox{d}u}
\newcommand{\dnu}{\hbox{d}\nu}
\newcommand{\mb}[1]{\mathbb{#1}}
\newcommand{\prf}{\smallskip \noindent{\it Proof}.\ }

\title{Solution to the Inverse Wulff Problem by Means of the Enhanced Semidefinite Relaxation Method}

\headlinetitle{Solution to the inverse Wulff problem}

\lastnameone{\v{S}ev\v{c}ovi\v{c}}
\firstnameone{Daniel}
\nameshortone{D. \v{S}ev\v{c}ovi\v{c}}
\addressone{Dept. of Applied Mathematics and Statistics, Comenius University, Mlynsk\'a dolina, 842 48 Bratislava}
\countryone{Slovakia}
\emailone{sevcovic@fmph.uniba.sk}

\lastnametwo{Trnovsk\'a}
\firstnametwo{M\'aria}
\nameshorttwo{M. Trnovsk\'a}
\addresstwo{Dept. of Applied Mathematics and Statistics, Comenius University, Mlynsk\'a dolina, 842 48 Bratislava}
\countrytwo{Slovakia}
\emailtwo{trnovska@fmph.uniba.sk}

\abstract{
We propose a novel method of resolving the optimal anisotropy function. The idea is to construct the optimal anisotropy function as a 
solution to the inverse Wulff problem, i.~e. as a minimizer for the anisoperimetric ratio for a given Jordan curve in the plane. It leads to a nonconvex quadratic optimization problem with linear matrix inequalities. In order to solve it  we propose the so-called enhanced semidefinite relaxation method which is based on a solution to a convex semidefinite problem obtained by a semidefinite relaxation of the original problem augmented by quadratic-linear constraints. We show that the sequence of finite dimensional approximations of the optimal anisoperimetric ratio converges to the optimal anisoperimetric ratio which is a solution to the inverse Wulff problem. Several computational examples, including those corresponding to boundaries of real snowflakes and discussion on the rate of  convergence of numerical method are also presented in this paper.}



\researchsupported{This research was supported by the projects VEGA 1/2429/12 and 7FP EU STRIKE No. 304617}

\begin{document}

\section{Introduction}

The classical isoperimetric inequality $L^2 \ge 4\pi {\mathcal A}$ relates the length $L$ of a Jordan curve $\Gamma$ in the plane $\R^2$ and the area ${\mathcal A}$ enclosed by $\Gamma$. The equality is attained if and only if  $\Gamma$ is a circle. It was apparently known to antique mathematician  Pappus from Alexandria (c.f. \cite{federico}). In \cite{Wulff1901} Wulff formulated and later Dingas in \cite{Dinghas1944} rigorously proved the isoperimetric inequality in the framework of the so-called relative Finsler geometry. Given the Jordan curve $\Gamma$ in the plane and the Finsler metric function $\Phi$ we can define the total interface energy functional $L_{\Phi}(\Gamma)=\int_{\Gamma}\Phi(\vecn)\,\ds$ where $\vecn$ is the unit inward normal vector to $\Gamma$. Then an analogous anisoperimetric inequality is satisfied for the so-called anisoperimetric ratio $\Pi_\Phi$ (or the isoperimetric ratio in the relative Finsler geometry given by the metric $\Phi$)
\[
\Pi_\Phi(\Gamma) := \frac{L_{ \Phi}(\Gamma)^2}{4|W_{ \Phi}| {\mathcal A}(\Gamma)} \ge 1.
\]
Here $|W_\Phi|$ is the area of the Wulff shape corresponding to the Finsler metric $\Phi$.  The anisoperimetric inequality has been proved and generalized to any spatial dimension (see \cite{dacorogna}).

Knowledge of the Finsler metric function $\Phi$ plays an essential role in many applied problems. In particular, in material science the Finsler metric function enters many crystal growth models based on Allen-Cahn type of nonlinear parabolic partial differential equations (c.f. \cite{benes2003,deckelnik2002,dziuk1999,hong2013,pozzi2008} and other references therein). In \cite{belletini} Bellettini and Paolini derived the Allen-Cahn parabolic partial differential equation for the gradient flow for the anisotropic Ginzburg-Landau free energy 
\[
{\mathcal E}(u) = \int_\Omega \frac{\xi}{2}\Phi(\nabla u)^2 +\frac{1}{\xi}f(u) \hbox{d}x,
\]
where $\Phi$ is the Finsler metric function. Here the function $u\in[-1,1]$ stands for the order parameter characterizing two phases ($u=\pm1$) of a material. The function $f$ is a double-well potential that gives rise to a phase separation and $\xi \ll 1$ is a small parameter representing thickness of the interface (c.f. \cite{graser}). Another important application involving the anisotropy function arises from motion of planar interfaces in which a family of curves is evolved in the normal direction by the velocity
\[
v = \kappa_{\Gamma,\Phi} + f,
\]
where $\kappa_{\Gamma,\Phi}$ is the so-called anisotropic curvature (c.f. \cite{dziuk1999,deckelnik2002,belletini,benes2003,hong2013} and Section~2.2). Such a flow also has a special importance in anisotropic diffusion image segmentation and edge detection models (see \cite{perona,weickert,strachota,MS2004}). Knowing underlying image anisotropy one can construct an efficient algorithm to segment important boundaries in the image or even denoising it by means of a anisotropic variant of Perona-Malik model \cite{perona,weickert}. 

However, less attention is put on understanding and resolving the Finsler ani\-sotropy function itself. The main purpose of this paper is to propose a novel method of determining the optimal Finsler metric function. The main idea is to resolve the Finsler metric with respect to a given planar curve representing thus a benchmark for underlying anisotropy. For instance, a boundary of a snowflake can be considered as such a benchmark curve yielding the optimal Finsler metric for its crystal growth model. In our approach, the idea is to find the underlying anisotropy function $\Phi$ by means of minimization of the anisoperimetric ratio. Due to properties of anisoperimetric ratio this approach can be viewed as a method of construction of the Finsler metric $\Phi$ that minimizes the total interface energy $L_\Phi$ for a given Jordan curve $\Gamma$ in the plane provided that the area of the Wulff shape is prescribed. It can be also regarded as a solution to the inverse Wulff problem stated as follows: given a Jordan curve $\Gamma$, find an optimal anisotropy function $\Phi$ minimizing the anisoperimetric ratio, i.~e.
\[
\inf_\Phi \Pi_\Phi(\Gamma).
\]
In this paper we show how to solve the inverse Wulff problem by means of nonconvex optimization and semidefinite relaxation methods and techniques. We will reformulate the inverse Wulff problem in terms of an indefinite quadratic optimization problem with linear matrix inequality constraints. It is shown that this problem belongs to a general class of quadratic optimization problems with linear and semidefinite constraints. In the proposed method of enhanced semidefinite relaxation, an equality constraint of the form $Ax=b$ are augmented by the quadratic-linear constraint $Ax x^T = b x^T$. Although it is a dependent constraint, it turns out that semidefinite relaxation of such an augmented problem leads to a convex semidefinite program (SDP) obtained as a second Lagrangian dual problem to the augmented indefinite quadratic optimization problem. Since the convexity of SDP is enhanced by the augmented quadratic-linear constraint we will refer to this method as the enhanced semidefinite relaxation method. The resulting SDP can be efficiently solved by using of available solvers for nonlinear programming problems over symmetric cones, e.g. SeDuMi or SDPT3 Matlab solvers \cite{sturm}. The method of the enhanced semidefinite relaxation can be also used in other applications leading to nonconvex constrained problems. 

The paper is organized as follows. In the next section we introduce necessary notation. 
We also recall known facts from parametric description of planar curves, Finsler relative geometry and anisoperimetric inequality. Section~3 is devoted to the Fourier series representation of the inverse problem. We also provide two useful criteria for nonnegativity of the Fourier series expansion given in terms of positive semidefinite Toeplitz matrices. In Section~4 we introduce and investigate properties of the Fourier length spectrum of a planar curve. We investigate its useful properties and derive important estimates. We furthermore reformulate the optimization problem in terms of Fourier coefficients of the anisotropy function. Section~5 is devoted to a method of the enhanced semidefinite relaxation of nonconvex quadratic optimization problem. We derive relatively simple sufficient conditions under which the primal problem and its semidefinite relaxed problem yield the same optimal value. Analysis of convergence of finite dimensional approximations is studied in Section~6. Finally, in Section~7 we present several computational examples illustrating optimal anisotropy functions minimizing the anisoperimetric ratio for various classes of planar curves including in particular examples of snowflakes. We also  investigate the experimental order of complexity and convergence of finite dimensional approximations to the solution of the inverse Wulff problem. 

\section{Preliminaries and notations}

\subsection{Parameterization of plane curves}

Following \cite{MS2001,SY2012} we introduce a notation for parameterization of planar curves. 
Let $\Gamma\subset \R^2$ be a $C^1$ smooth curve of a finite length, i.~e. $\Gamma$ can be parameterized by a $C^1$ mapping 
$\vecx:[0,1]\to \R^2$, $\Gamma=\{ \vecx(u), u\in [0,1]\}$ such that $\Vert\partial_u \vecx(u)\Vert >0$ for any $u\in [0,1]$. 
Here $\Vert a\Vert=\sqrt{a^T a}$ denotes the Euclidean norm of a vector $a$. For a $C^1$ smooth Jordan curve $\Gamma$ 
(simple and closed curve in the plane) we will assume the parameterization $\vecx$ of $\Gamma$ is counterclockwise and we will 
impose periodic boundary conditions for $\vecx(u)$ at $u=0,1$. For each point $\vecx=\vecx(u)\in\Gamma$ we can define the unit tangent 
vector $\vect=\partial_u \vecx/\Vert\partial_u \vecx\Vert$. The tangent angle $\nu$ can be defined through the relation 
$\vect\equiv(t_1,t_2)^T=(\cos\nu, \sin\nu)^T$. The unit inward vector $\vecn$ then satisfies $\vecn=(-\sin\nu, \cos\nu)^T$. 
The arc-length parameterization $s\in[0,L(\Gamma)]$ is related to the fixed domain parameterization $u\in[0, 1]$ by the relation: 
$\ds = \Vert\partial_u \vecx\Vert \du$. Then $\vect=\partial_s\vecx$. For a $C^2$ smooth curve $\Gamma$ we also recall the Frenet formulae: $\partial_s\vect = \kappa\vecn$ and $\partial_s\vecn = -\kappa\vect$ where $\kappa=\det(\partial_s\vecx, \partial_s^2\vecx)$ 
is the curvature. For the tangent angle $\nu$ we obtain $\partial_s\nu =\kappa$. Notice that  for a strictly convex curve $\Gamma$ 
(i.~e. $\hbox{int}(\Gamma)$ is strictly convex) the sign of the curvature $\kappa$ is positive and so the tangent angle $\nu\in[0,2\pi]$ 
can be used as a parameterization of $\Gamma$ and $\dnu =\kappa \ds$.
The total length $L(\Gamma)$ and the area ${\mathcal A}(\Gamma)$ enclosed by a $C^1$ Jordan curve $\Gamma$ are given by 
$L(\Gamma) = \int_\Gamma \ds$ and ${\mathcal A}(\Gamma) = -\frac{1}{2}\int_\Gamma\vecx^T \vecn\,\ds$.

\subsection{Finsler metric and description of the relative geometry}

In this section we recall basic facts and notations regarding description of the relative Finsler geometry in $\R^n$. Following Paolini \cite{paolini1998} and Gr\"aser \cite[Assumptions A1,A2]{graser}, we consider the so-called Finsler metric $\Phi:\R^n\to\R_+$ has the following properties:
\begin{enumerate}
\item $\Phi$ is a positively homogeneous function of degree one, \\ i.~e. $\Phi(t\vecx) = t\Phi(\vecx)$ for each $\vecx\in\R^n$ and $t\ge0$;
\item $\Phi$ is a $C^2$ smooth function and $\Phi(x)>0$ in $\R^n\setminus\{0\}$, $\Phi(0)=0$;
\item $\Phi^2$ is a strictly convex function.
\end{enumerate}
Regularity assumptions on the Finsler metric $\Phi$ have been discussed in Bene\v{s} et al. \cite{benes2004}.

\begin{figure}
\begin{center}
\subfloat[]{\includegraphics[width=0.28\textwidth]{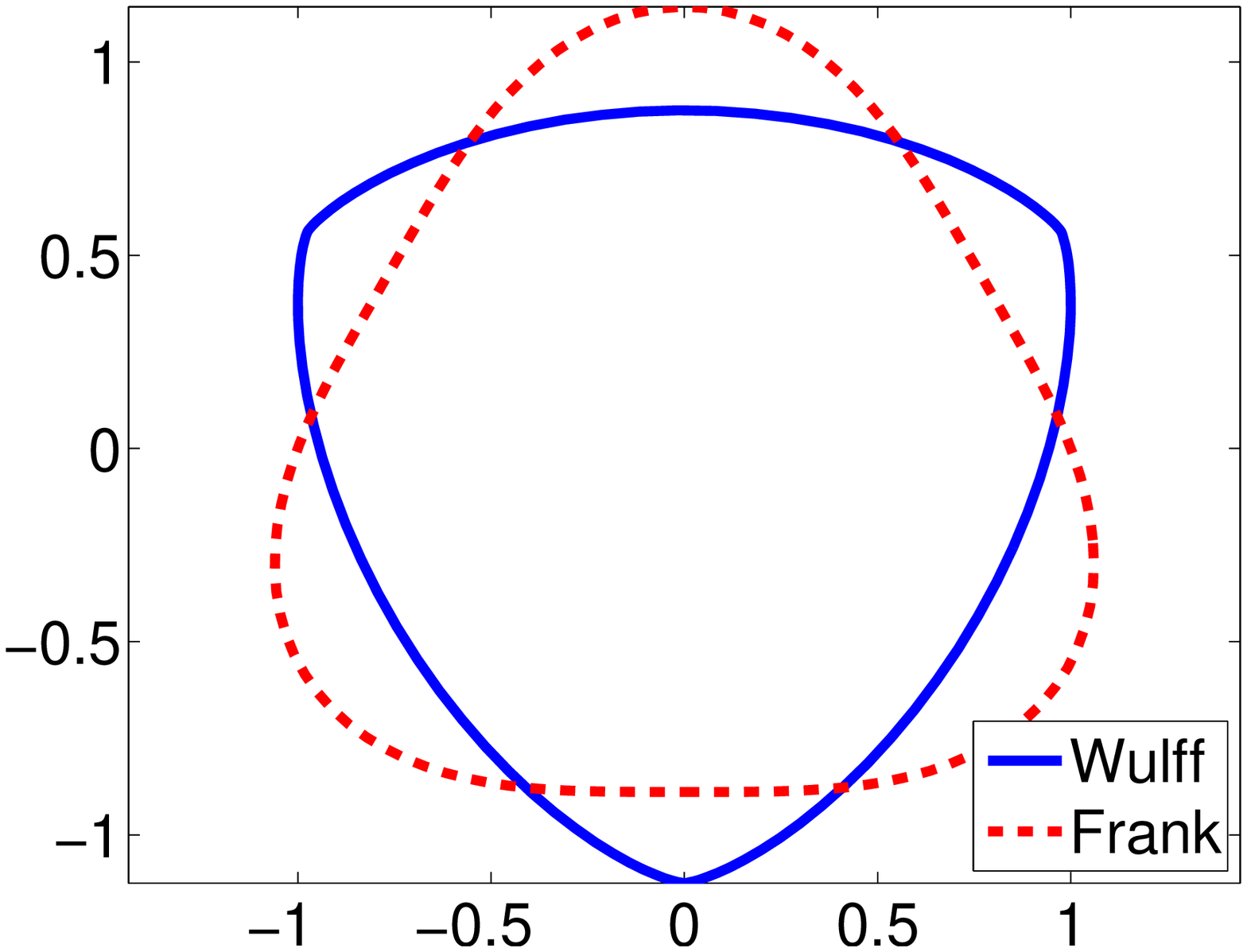}}
\subfloat[]{\includegraphics[width=0.28\textwidth]{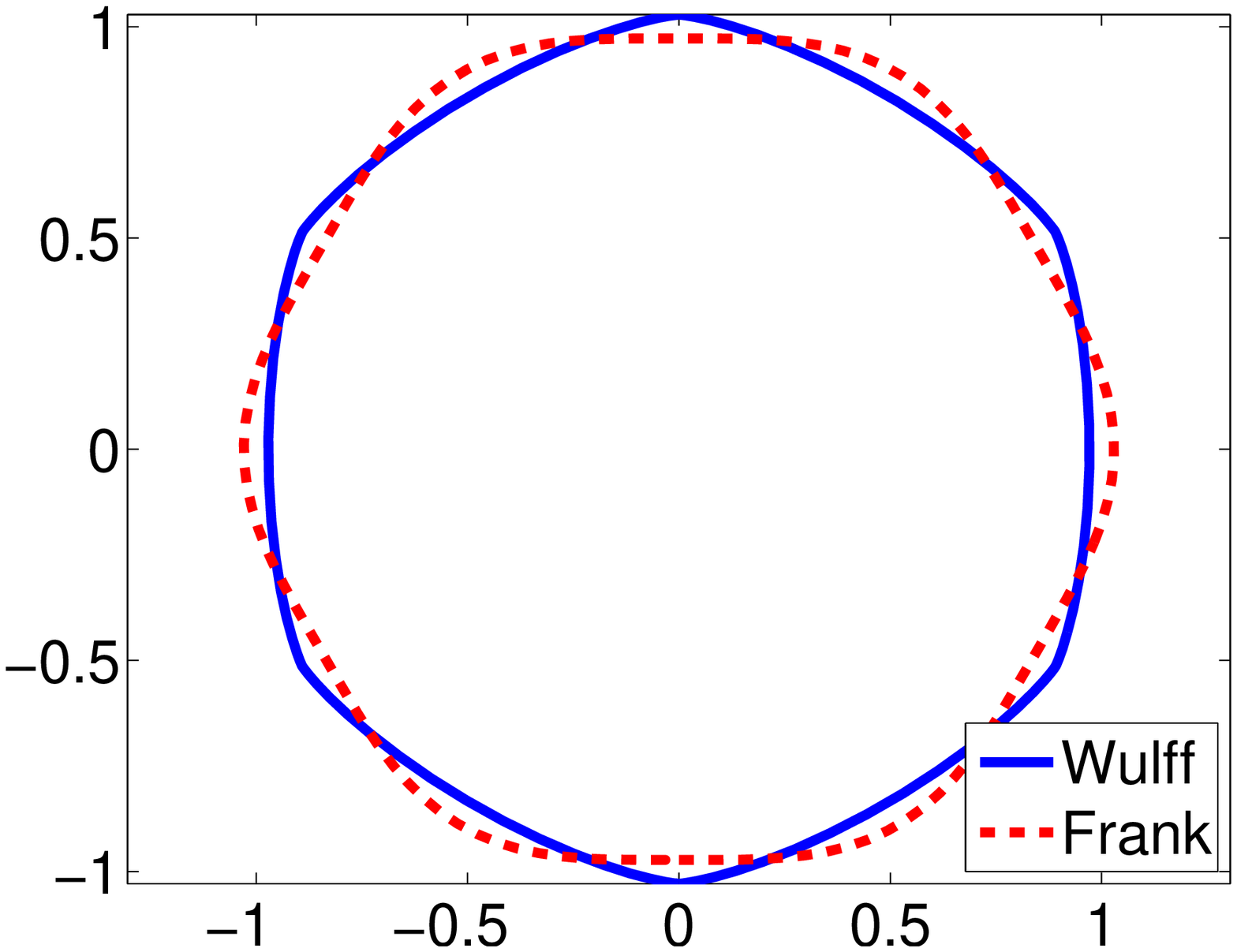}}
\subfloat[]{\includegraphics[width=0.28\textwidth]{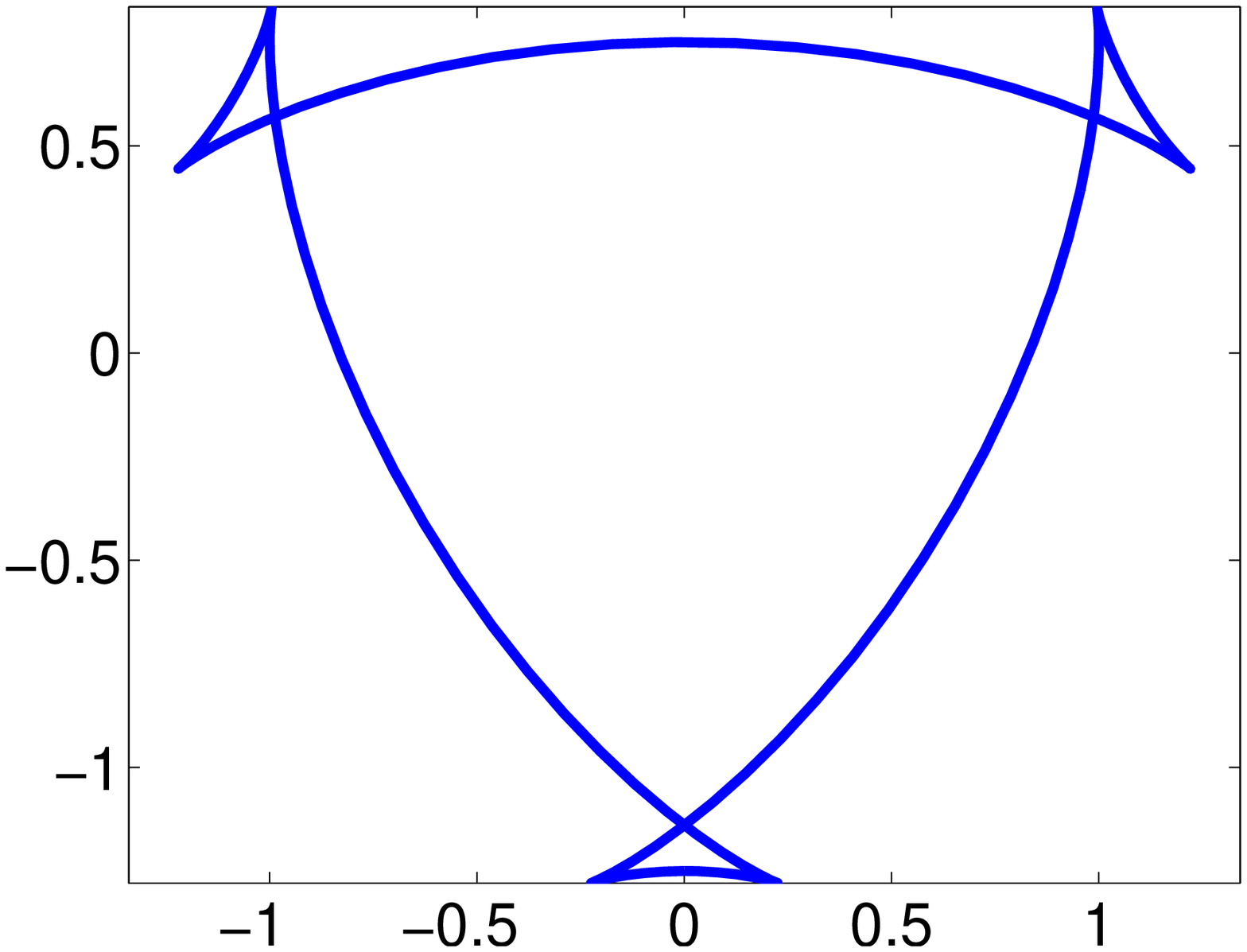}}
\end{center}
\caption{\small
(a) the Wulff shape and Frank diagram for the anisotropy function $\sigma(\nu)=1+\varepsilon\cos(m\nu)$ for the three-fold anisotropy $m=3$; (b) Hexagonal anisotropy $m=6$, $\varepsilon=1/(m^2-1)$; (c) a parametric curve \eqref{boundaryWulff} with tails for $m=3, \varepsilon=1/4$.
}
\label{fig:example}
\end{figure}

\begin{rem}
In classical definitions of the Finsler metric (c.f. \cite{belletini,dziuk1999}), absolute homogeneity property of $\Phi$, i.~e. $\Phi(t\vecx) = |t|\Phi(\vecx)$ for each $\vecx\in\R^n$ and $t\in\R$, is usually assumed. In contrast to such an assumption on absolute homogeneity of $\Phi$, in our definition we allow $\Phi$ to belong to a larger class of functions. In particular, we consider a class of anisotropy functions having odd number of folds  (c.f. \cite{chalmers}, see also \cite{kobayashi,graser}). For example, a three-fold anisotropy function depicted in Figure~\ref{fig:example}, can be found as a shape of the (111) facet of Pb particles, prepared and equilibrated on Cu(111) under ultrahigh vacuum conditions (c.f. \cite{surnev})
\end{rem}

Then the Wulff shape $W_\Phi$ and the Frank diagram ${\mathcal F}_\Phi$ corresponding to the Finsler metric $\Phi$ can be defined as follows:

\begin{equation}
W_\Phi = \bigcap_{\Vert\vecn\Vert=1}\{\vecx \in \R^n\ |\  -\vecx^T\vecn \le \Phi(\vecn)  \},\ \ {\mathcal F}_\Phi = \{\vecx \in \R^n\ |\  \Phi(-\vecx)\le 1  \}.
\label{wulff-grank-phi}
\end{equation}

The Wulff shape $W_\Phi$ is always a convex set. In the case $\Phi(\vecx)=\Vert\vecx\Vert$ the Wulff shape and the Frank diagram are just unit balls in $\R^n$. 
If we restrict our attention to the plane $\R^2, n=2,$ we can provide a simplified  characterization of the Finsler metric $\Phi$ by means of the real anisotropy function $\sigma=\sigma_\Phi$ where
\begin{equation}
\sigma_\Phi(\nu) = \Phi(\vecn), \quad \hbox{where} \ \vecn=(-\sin\nu, \cos\nu)^T,
\label{sigma-phi}
\end{equation}
and vice versa, with regard to \eqref{wulff-grank-phi}, the Finsler metric $\Phi$ can be constructed from $\sigma$ as follows: 
\[
\Phi_\sigma(-\vecx) = \sigma(\nu)  \Vert\vecx\Vert, \quad\hbox{where}\ \vecx/\Vert\vecx\Vert = -\vecn, \ \ \hbox{for}\ \vecx\not= 0, \Phi_\sigma(0)=0.\] 
The anisotropy function $\sigma:\R\to\R^+_0$ is assumed to be a $2\pi$-periodic smooth function of the tangent angle $\nu$. If we restrict our attention to the class of $\pi$-periodic anisotropy functions then the corresponding Finsler metric $\Phi_\sigma$ is an absolute homogeneous function. 

In terms of the anisotropy function $\sigma$, the Wulff shape $W_\sigma$ and the Frank diagram can be described as follows:
\[
W_{\sigma}=\bigcap_{\nu\in[0,2\pi]}\left\{\vecx\in\R^2 \ |  \ -\vecx^T \vecn \leq \sigma(\nu)\right\},
\ \ 
{\mathcal F}_\sigma = \{\vecx= - r \vecn \ |\  0\le r\le 1/\sigma(\nu)  \}.
\]
where $\vect=(\cos(\nu), \sin(\nu))^T$ and $\vecn=(-\sin(\nu), \cos(\nu))^T$ are unit tangent and inward normal vectors. Since $\partial_\nu\vect = \vecn$ and $\partial_\nu\vecn = - \vect$, its boundary $\partial W_{\sigma}$ can be parameterized as
\begin{equation}
\partial W_{\sigma}=\left\{\vecx(\nu)\ |\ \vecx(\nu)=-\sigma(\nu)\vecn+\sigma'(\nu)\vect, \ \nu\in  [0,2\pi]\right\}
\label{boundaryWulff}
\end{equation}
provided that the Frank diagram ${\mathcal F}_\sigma$ is strictly convex, i.~e. $\Phi_\sigma$ is strictly convex and smooth. Notice that the right hand side of \eqref{boundaryWulff} is the set of all Cahn-Hoffman vectors of the form $\vecx=-\nabla \Phi_\sigma(\vecn), \Vert\vecn\Vert=1$.

As $\dnu=\kappa\ds$ we have $\vect=\partial_s\vecx = \kappa\partial_\nu\vecx$ and $\partial_s^2\vecx=\partial_s\vect=\kappa\vecn$. Since $\partial_\nu\vecx = (\sigma+\sigma'')\vect - 2\sigma'\vecn$ for $\partial W_{\sigma}$ we obtain $\kappa=\det(\partial_s\vecx,\partial_s^2\vecx)= \kappa^2 (\sigma+\sigma'')$ and so the curvature $\kappa$ of $\partial W_{\sigma}$ is given by $\kappa=[\sigma(\nu)+\sigma''(\nu)]^{-1}$ (c.f. \cite{SY}). Hence the Wulff shape $W_{\sigma}$  is strictly convex if and only if $\sigma+\sigma'' > 0$. If we define the anisotropic curvature by the relation: ${\kappa_\sigma}:= [\sigma(\nu)+\sigma''(\nu)] \kappa$ then $\kappa_\sigma\equiv 1$ on $\partial W_{\sigma}$. 

Finally, the area $|W_\sigma|$ of the Wulff shape entering the anisoperimetric ratio can be calculated as follows:
\begin{eqnarray}
|W_\sigma| &=& -\frac{1}{2}\int_{\partial W_{\sigma}}\vecx^T \vecn\,\ds
=\frac{1}{2}\int_{\partial W_{\sigma}}\sigma(\nu)\,\ds
\label{areaW} \\
&= &\frac12 \int_0^{2\pi} \sigma(\nu) [\sigma(\nu) + \sigma''(\nu)] \hbox{d}\nu = \frac12 \int_0^{2\pi} |\sigma(\nu)|^2 - |\sigma'(\nu)|^2 \hbox{d}\nu
\nonumber
\end{eqnarray}
because $\dnu = \kappa \ds = [\sigma + \sigma'']^{-1} \ds$. Clearly, if $\sigma\equiv 1$ then the boundary $\partial W_{1}$ of $W_1$ is a circle with the radius 1, and $|W_{1}|=\pi$.

In Figure~\ref{fig:example} we show typical examples of the anisotropy functions with three-fold and hexagonal symmetries. We consider a class of anisotropy functions of the form $\sigma(\nu)=1+\varepsilon\cos(m\nu)$ where $m\in\N$ (c.f. \cite{chalmers,kobayashi}). The parameter $\varepsilon\ge 0$ represents the so-called strength of anisotropy. Clearly, $\sigma\ge 0, \sigma+\sigma''\ge 0$ provided that $\varepsilon\le 1/(m^2-1)$. In Figure~\ref{fig:example} (c) we plot the curve $\partial W_\sigma$ given by \eqref{boundaryWulff} for the case $m=3$ and $\varepsilon=1/4>1/(m^2-1)$. In such a case $\sigma(\nu)+\sigma''(\nu)$ becomes negative for some angles $\nu$ and $\partial W_\sigma$ is no longer a Jordan curve. Hence the condition $\sigma+\sigma''\ge 0$ is crucial for the analysis of a Wulff shape.

\subsection{Anisoperimetric ratio and anisoperimetric inequality}

In \cite{SY} Yazaki and the first author proved the mixed anisoperimetric inequality:
\begin{equation}
\frac{L_\sigma(\Gamma) L_\mu(\Gamma)}{{\mathcal A}(\Gamma)} \ge K_{\sigma,\mu}, \qquad\hbox{where}\  K_{\sigma,\mu} = 2 \sqrt{|W_\sigma| |W_\mu| } + L_\sigma (\partial W_\mu).
\label{mixedanisoperim}
\end{equation}
It has been shown for a $C^2$ smooth Jordan curve $\Gamma$ and an arbitrary pair of anisotropy functions $\sigma, \mu\in\K$ belonging to the cone of $2\pi$-periodic functions that 
\begin{equation}
\K=\{\sigma\in W^{2,2}_{per}(0,2\pi)\ |\ \sigma(\nu)\ge 0, \sigma(\nu) + \sigma''(\nu) \ge 0, \  \hbox{for a.e.}\  \nu\in[0,2\pi]\}.
\label{cone}
\end{equation} 
The minimum in inequality \eqref{mixedanisoperim} is attained for a curve $\Gamma$ which is a certain convex combination 
of boundaries $\partial W_\sigma$ and $\partial W_\mu$ of Wulff shapes (c.f. \cite[Theorem 2]{SY}). Here  $W^{r,2}_{per}(0,2\pi)$ 
denotes the Sobolev space of all real valued $2\pi$-periodic functions having their distributional derivatives square integrable up to 
the order $r$. This is a Hilbert space when endowed by the norm 
$\Vert\sigma\Vert_{r,2}= \left(\sum_{k=0}^\infty (1+ k^{2r}) |\sigma_k|^{2}\right)^\frac12$. By $\{ \sigma_k\in\CC, k\in \mathbb{Z}\}$ 
we have denoted the set of complex Fourier coefficients of the function $\sigma$ (see \eqref{fourier} below). In the particular 
case when $\mu=\sigma\in\K$ we have $K_{\sigma,\sigma} = 2 |W_\sigma| + L_\sigma (\partial W_\sigma) = 4 |W_\sigma|$ because 
$ L_\sigma (\partial W_\sigma)=\int_{\partial W_{\sigma}}\sigma\,\ds =  2|W_\sigma|$ (see \eqref{areaW}). Therefore, as a consequence 
of \eqref{mixedanisoperim} we obtain the anisoperimetric inequality 
\begin{equation}
\Pi_\sigma(\Gamma) \equiv \frac{L_{\sigma}(\Gamma)^2}{4|W_{\sigma}| {\mathcal A}(\Gamma)} \ge 1.
\label{anisoineq}
\end{equation}
The equality is attained if and only if $\Gamma$ is homothetically similar to $\partial W_\sigma$. 
Note that both \eqref{mixedanisoperim} and its special case \eqref{anisoineq} are generalizations of the anisoperimetric 
inequality by Wulff \cite{Wulff1901} (see also \cite{dacorogna}) to the case of $2\pi$-periodic anisotropy functions, 
i.~e. for a larger class of positive homogeneous Finsler metric function. It is worth noting that the area $|W_\sigma|$ 
of a Wulff shape satisfies:
\[
0\le | W_\sigma| = \frac{L_\sigma(\Gamma)^2}{4\Pi_\sigma(\Gamma){\mathcal A}(\Gamma)}\le \frac{L_\sigma(\Gamma)^2}{4 {\mathcal A}(\Gamma)} <\infty
\]
for any anisotropy function $\sigma\in\K$ and a Jordan curve $\Gamma$.

For a given smooth Jordan curve $\Gamma$ in the plane, our goal is to find the optimal anisotropy function $\sigma$ minimizing 
the anisoperimetric ratio:
\begin{equation}
\inf_{ \sigma\in\K} \Pi_\sigma(\Gamma).
\label{anisoproblem}
\end{equation}
It is useful to emphasize that the following homogeneity conditions hold true:
\begin{equation}
L_{t\sigma}(\Gamma) = t L_{\sigma}(\Gamma), \quad |W_{t\sigma}| = t^2 |W_\sigma|, \quad 
\Pi_{t\sigma}(\Gamma) = \Pi_{\sigma}(\Gamma), 
\label{homog}
\end{equation}
for any $\sigma\in \K$ and all  $t>0$. The anisoperimetric ratio is therefore a homogeneous function of the zero-th order with 
respect to positive scalar multiple of the anisotropy function $\sigma$. In order to solve problem (\ref{anisoproblem}) uniquely 
with respect to scalar multiples of $\sigma$, instead of (\ref{anisoproblem}), we can solve the maximization problem
\begin{equation}\label{maxW}
\begin{array}{rl}
\sup & |W_\sigma| \\
{\rm s. t.} & L_\sigma(\Gamma)=L(\Gamma),
\\
& \sigma \in \K.
\end{array} 
\end{equation}
It means that the goal is to maximize the area of the Wulff shape under the constraint that the interface energy $L_\sigma(\Gamma)$ 
is fixed to the constant length $L(\Gamma)$. The choice of the scaling constraint $L_\sigma(\Gamma)=L(\Gamma)$ is quite natural 
because in the case $\Gamma$ is a circle, the anisotropy function $\sigma\in \K$ maximizing $|W_\sigma|$ under the constraint  
$L_\sigma(\Gamma)=L(\Gamma)$ is just unity, $\sigma\equiv1$. It should be also obvious that, up to a positive multiple of $\sigma$, 
a solution $\sigma$ to the maximization problem \eqref{maxW} is also a solution to the minimal interface energy problem 
\begin{equation}\label{minL}
\begin{array}{rl}
\inf & L_\sigma(\Gamma) \\
{\rm s. t.} & |W_\sigma| = 1,
\\
& \sigma \in \K. 
\end{array} 
\end{equation}
Hence the problem of resolving the minimizer $\sigma$ for the anisoperimetric ratio can be also viewed as a problem of finding 
the anisotropy function minimizing the total interface energy.  

\section{Fourier series representation}

Let $\sigma:\R\to\R$ be a $2\pi$-periodic function, $\sigma\in W^{2,2}_{per}(0,2\pi)$. It can be represented by its complex Fourier series 
\begin{equation}
\sigma(\nu) = \sum_{k=-\infty}^\infty \sigma_k e^{i k\nu}, \quad \hbox{where} \ 
\sigma_k=\frac{1}{2\pi}\int_0^{2\pi} e^{-i k\nu} \sigma(\nu)\hbox{d}\nu
\label{fourier}
\end{equation}
are complex Fourier coefficients. Since $\sigma(\nu)$ is assumed to be a real function we have $\sigma_{-k} = \bar\sigma_k$ for any $k\in\Z$ and $\sigma_0\in\R$. Notice that for $\sigma\in W^{2,2}_{per}(0,2\pi)$ we have $\sigma(\nu)+\sigma''(\nu) = \sum_{k=-\infty}^\infty (1-k^2)\sigma_k e^{i k\nu}$ in the norm of the Lebesgue space $L^{2}(0,2\pi) = W^{0,2}_{per}(0,2\pi)$. 

It follows from \eqref{areaW} that the area $|W_\sigma|$ of the Wulff shape can be  expressed in terms of Fourier coefficients as follows:
\begin{eqnarray}
|W_\sigma| &=& \frac12 \int_0^{2\pi} |\sigma|^2 - |\sigma'|^2 \hbox{d}\nu 
= \frac12 \int_0^{2\pi} \sum_{k,m=-\infty}^\infty \bar\sigma_m \sigma_k (1- m k) e^{i(k-m)\nu} \dnu
\nonumber\\
&=& \pi \sum_{k=-\infty}^\infty (1-k^2) |\sigma_k|^2 = \pi\sigma_0^2 + 2\pi \sum_{k=1}^\infty (1-k^2) |\sigma_k|^2. 
\label{area-Wulff}
\end{eqnarray}

Similarly, we can express the interface energy 
\begin{equation}
L_\sigma(\Gamma) = \int_\Gamma \sigma(\nu)\ds = \sum_{k=-\infty}^\infty \sigma_k \int_\Gamma e^{i k \nu}\ds
= c_0 \sigma_0 + 2\Re  \sum_{k=1}^\infty \bar c_k \sigma_k, 
\label{interface-energy}
\end{equation}
where the complex coefficients $c_k=\int_\Gamma e^{-i k \nu}\ds, k\in\Z$, form the so-called Fourier length spectrum of the curve $\Gamma$ (see Section 4).

\subsection{Criteria for nonnegativity of Fourier series}
In this section we recall two useful criteria guaranteeing nonnegativity of complex Fourier series. Both of them are based on positive definiteness of certain Hermitian matrices related to the complex Fourier coefficients. 

For a transpose of the matrix $A$ we will henceforth write $A^T$. For a complex conjugate of a complex matrix $H$ we will write $H^*$, i.~e. $H^*=\bar{H}^T$. The sets of real $N\times N$ symmetric and complex Hermitian matrices are denoted by $\mathcal{S}^N$ and $\mathcal{H}^N$, respectively, i.~e. $\mathcal{S}^N=\{ A\in \mb{R}^{N\times N} \ | \ A=A^T\}, \quad \mathcal{H}^N=\{ H\in \mb{C}^{N\times N} \ | \ H=H^*\}$.
We will write $A\succeq 0$ ($A\succ 0$), if a real symmetric matrix or a complex  Hermitian matrix $A$ is positive semidefinite (positive definite). That is, $A\succeq 0$ ($A\succ 0$) if $x^TAx\geq 0$ 
($x^TAx>0$) for all $x\in \mathbb{R}^N, x\neq 0$, if the real case and if
$z^*Az\geq 0$ ($z^*Az>0$) for all $z\in \mathbb{C}^N, z\neq 0$, if the complex case.

An easy criterion is based on the Bochner theorem on positive-definite real functions. Given Fourier series representation (\ref{fourier}) of a 
smooth $2\pi$-periodic function $\sigma:\R\to\R$ we can construct the following Toeplitz circulant matrix:
\begin{equation}
Q=\hbox{Toep}(\sigma_0, \sigma_1, \ldots, \sigma_{N-1})=
\left(
\begin{matrix}
\sigma_0 & \bar\sigma_{1} & \ldots & \bar\sigma_{N-1} \\
\sigma_1 & \sigma_{0} & \ldots  & \bar\sigma_{N-2} \\
 \vdots& \vdots & \ddots  & \vdots \\
\sigma_{N-1} & \sigma_{N-2} & \ldots  & \sigma_{0}
\end{matrix}
\right),
\end{equation}
i.~e. $Q_{pq} = \sigma_{p-q}$, where $\sigma_{-k}=\bar\sigma_k$. The matrix $Q$ is Hermitian, $Q\in\mathcal{H}^N$. 

\begin{tvr}\label{Dumit-pos} \cite[Th. 1.8]{dumit}
Let $\sigma:\R\to\R$ be a smooth $2\pi$-periodic function. Then $\sigma(\nu)\ge 0$ for any $\nu\in\R$ if and only if the Toeplitz matrix 
$Q^{(N)}=\hbox{Toep}(\sigma_0, \sigma_1, \ldots, \sigma_{N-1}) \succeq 0$ is positive semidefinite for any $N\in\N$. 
\end{tvr}

This is a relatively simple criterion. Its proof is rather straightforward and it is based on a similar argument as the one of Prop.~\ref{Fspec-prop}. Unfortunately, it includes infinitely many conditions $Q^{(N)}\succeq 0$ for any $N\in\N$ even in the case Fourier series expansion (\ref{fourier}) is finite. Indeed, let us consider the function $\sigma(\nu) = 1-a\cos(\nu)$. Clearly,  $\sigma_0=1, \sigma_{\pm1}= -a/2, \sigma_{\pm k}=0$ for $k\ge 2$. Then for any finite $N$ there exists $a>1$ such that  
$\hbox{Toep}(1, -a/2, 0, \ldots, 0)\in  \mathcal{H}^N$ is a positive semidefinite matrix but the function $\sigma(\nu)$ attains negative values. On the other hand, when $N\to+\infty$ then the range for $a$ is being restricted to $[-1,1]$, i.~e. $\sigma(\nu)\ge 0$.

However, the condition $\sigma(\nu)+\sigma''(\nu)\ge 0$ is crucial for avoiding of selfintersection of the parametric description \eqref{boundaryWulff} of the boundary $\partial W_\sigma$ of the Wulff shape. In \cite{mclean} McLean derived another useful criterion for nonnegativity of a partial finite  Fourier series sum. It is again formulated in terms of positive semidefinite Hermitian matrices. This criterion is a consequence of the classical Riesz-Fejer factorization theorem (c.f. \cite[pp. 117--118]{riesz}) and it reads as follows:

\begin{tvr}\label{Lean-pos} \cite[Prop. 2.3]{mclean}
Let $\sigma_0\in\R, \sigma_k=\bar\sigma_{-k}\in\CC$ for $k=1, \ldots, N-1$. Then the finite Fourier series expansion  $\sigma(\nu) = \sum_{k=-N+1}^{N-1} \sigma_k e^{i k\nu}$ is a nonnegative function $\sigma(\nu)\ge 0$ for $\nu\in\R$, if and only if the set $\mathcal{F}\subset  \mathcal{H}^N$ is nonempty, where 
\[
\mathcal{F}=\{ F\in \mathcal{H}^N \ |\ F\succeq 0,\ \sum_{p=k+1}^N F_{p,p-k} = \sigma_k, \ \hbox{for each}\ k=0, 1,\ldots, N-1\}.
\]
\end{tvr}

\subsection{Reformulation as a nonconvex quadratic optimization problem with semidefinite constraints}

In order to compute the optimal anisotropy function $\sigma$ as a limit of its finite Fourier modes approximation we introduce the finite dimensional subcone $\K^N$ of $\K$ where
\begin{equation}
\K^N =\{ \sigma \in \K\ |\  \ \exists (\sigma_0, \sigma_1, \ldots,\sigma_{N-1})^T\in\CC^N,\ \ 
\sigma(\nu) = \sum_{k=-N+1}^{N-1} \sigma_k e^{i k \nu}\},
\label{coneN}
\end{equation}
where $\sigma_{-k}=\bar\sigma_k$. We will identify the cone $\K^N$ with a cone in $\CC^N$ consisting of all vectors $(\sigma_0, \sigma_1, \ldots,\sigma_{N-1})^T\in\CC^N$ representing  functions of the form $\sigma(\nu)=\sum_{k=-N+1}^{N-1} \sigma_k e^{i k \nu}\in\K$. Then, for $\sigma\in\K^N$ and $N\in\N$ , we have 
\[
|W_{\sigma}| = \pi\sigma_0^2 + 2\pi \sum_{n=1}^{N-1} (1-n^2) |\sigma_n|^2 , \qquad L_{\sigma}(\Gamma) = c_0 \sigma_0 + 2\Re  \sum_{n=1}^{N-1} \bar c_n \sigma_n.
\]
Given $N\in\N$, our purpose is to solve the following finite dimensional optimization problem
\begin{equation}\label{maxWN}
\begin{array}{rl}
\max & |W_{\sigma}| \\
{\rm s. t.} & L_{\sigma}(\Gamma)=L(\Gamma),
\\
& \sigma \in \K^N. 
\end{array} 
\end{equation}
An optimal solution to \eqref{maxWN} will be denoted by $\sigma^N$. It should be emphasize that the constraint $\sigma \in \K^N$ ensures $\sigma(\nu)+\sigma''(\nu)\ge 0$ for $\sigma=\sigma^N$. Without such a constraint 
the optimum of $\max_{\sigma} \{ |W_\sigma|\ |\  L_{\sigma}(\Gamma)=L(\Gamma), \ \sigma(\nu) = \sum_{k=-N+1}^{N-1} \sigma_k e^{i k \nu} \}$ may have kinks and selfintersections of the boundary $\partial W_\sigma$ of the optimal Wulff shape as it is shown in Figure~\ref{fig:example} (c).

By Prop.~\ref{Lean-pos} and taking into account that $\sigma(\nu)+\sigma''(\nu) = \sum_{k=-N+1}^{N-1} (1-k^2) \sigma_k e^{i k \nu}$ for any $\sigma\in\K^N$ we end up with the following representation of the cone $\K^N$:

\begin{lem}\label{coneKrepresentation}
We have $\sigma\in\K^N$ if and only if there exist $F, G\in \mathcal{H}^N, F,G \succeq 0,$ such that 
\[
\sum_{p=k+1}^N F_{p,p-k} = \sigma_k, \quad \sum_{p=k+1}^N G_{p,p-k} = (1-k^2) \sigma_k \quad \hbox{for any }\ k=0, \ldots, {N-1}.
\]
\end{lem}

Finally, we will rewrite problem  \eqref{maxWN} in terms of real and imaginary parts $x^R, x^I\in\R^N$ of a solution vector $\sigma\in\K^N$. For this purpose, we decompose $\sigma$ into its real and imaginary parts:
\[
x = \left(\begin{array}{c}
x^R \\ x^I
\end{array} \right)\equiv [x^R; x^I] \in \R^{n},\ \hbox{where}\  n=2N, \quad \hbox{and}\ \sigma_k = x^R_k + i x^I_k,
\]
and introduce the real $n\times n$ matrix $P_0$ as follows:
\begin{equation}
P_0 =\hbox{diag}(p_0, p_1, \ldots, p_{N-1}, q_0, q_1, \ldots, q_{N-1}), 
\label{P0def}
\end{equation}
where $p_0=q_0=-\pi, p_k=q_k= 2\pi (k^2-1)$ for $k\ge 1$. We also decompose the Fourier length spectrum $\{c_k, k\ge 0\}$, 
as follows: $\alpha_k  +i \beta_k = 2 c_k, k\ge1,  \ \alpha_0 = c_0, \beta_0=0$, where $\alpha=(\alpha_0, \ldots, \alpha_{N-1})^T, \beta=(\beta_0, \ldots, \beta_{N-1})^T \in \R^N$. Next we define a $2\times n$ real matrix $A$ and the vector $b\in\R^2$ as follows:
\begin{equation}
A=\left(
\begin{array}{cc}
 \alpha^T & \beta^T \\ 
 0_N^T & e_1^T
\end{array} 
\right), 
\qquad 
b= \left(
\begin{array}{c}
 L(\Gamma) \\ 
 0
\end{array} 
\right),
\label{Abdef}
\end{equation}
where $0_N=(0,\ldots,0)^T\in\R^N, e_1=(1,0,\ldots,0)^T\in\R^N$. With help of the semidefinite representation of the cone $\K^N$ from Lemma~\ref{coneKrepresentation} we can reformulate \eqref{maxWN} as follows:
\begin{equation}\label{maxWNP}
\begin{array}{rl}
\min & x^T P_0 x \\
{\rm s. t.} & A x = b, \qquad x \equiv [x^R; x^I],
\\ & 
\sum_{p=k+1}^N F_{p,p-k} = x^R_k + i x^I_k, \ \ \sum_{p=k+1}^N G_{p,p-k} = (1-k^2)(x^R_k + i x^I_k),
\\
& \hbox{for}\ \  k=0,\ldots, N-1,
\\ & F, G\succeq 0.
\end{array} 
\end{equation}
The equation from the second row in $Ax=b$  guarantees $x^I_0=0$, i.~e. $\sigma_0\in\R$. It is worth noting that the matrix $P_0$ is indefinite and this is why problem \eqref{maxWNP} is a nonconvex optimization problem with 
linear matrix inequality constraints. In Section~5 we will investigate a general class of nonconvex optimization problems of the form \eqref{maxWNP} and we will show that \eqref{maxWNP} can be solved by means of the enhanced semidefinite relaxation method based on the second Lagrangian dual to \eqref{maxWNP} augmented by a quadratic-linear constraint.

\section{The Fourier length spectrum of a curve}

In this section, we introduce a notion of the so-called complex Fourier length spectrum. It is related to Fourier series expansion of a quantity depending on the tangent angle $\nu$ of the unit tangent vector $\vect=(t_1,t_2)^T$. 

\begin{defi}\label{Fspec-def}
Let $\Gamma$ be a $C^1$ smooth curve in the plane. By the complex Fourier length spectrum of $\Gamma$ we mean the set $\{c_p, p\in\Z\}$ of all Fourier complex coefficients defined as follows:
\[
c_p = \int_\Gamma e^{-ip\nu}\ds = \int_\Gamma (t_1 - i t_2)^p \ds\,,
\]
where $\vect=(t_1,t_2)^T=(\cos(\nu), \sin(\nu))^T$ is the unit tangent vector to $\Gamma$.
\end{defi}

\medskip
\begin{ex}\label{fourier-spec-ex}
For example, if $\Gamma$  is a circle with a radius $r>0$ then we have $c_0=2\pi r$ and $c_k=0$ for each $k\not=0$. Indeed, $\Gamma$ can be parameterized by $x_1(u)=r\cos(2\pi u)$, $x_2(u)=r\sin(2\pi u)$. So $\nu=\pi/2 + 2\pi u$. Since $\ds = 2\pi r \du$ 
we have $c_k=0$ for each $k\not=0$. 
\end{ex}

\begin{ex}\label{fourier-spec-ex-capsule}
Let us consider a "capsule" curve consisting of two horizontal line segments with the length $l>0$ connected by half-arcs with a radius $r>0$ (see Figure~\ref{fig:capsule}). Then $c_0=L(\Gamma) = 2l + 2\pi r$ and $c_{2k}=2 l, c_{2k+1}=0$ because the tangent angle $\nu\in\{0,\pi\}$ on the line segments and integration of $e^{-i k\nu}$ over the union of remaining half-arcs yields zero for any $k\not=0$.
\end{ex}

\medskip
Concerning properties of the Fourier length spectrum we can formulate the following result.

\begin{tvr}\label{Fspec-prop}
Let $\Gamma$ be a $C^1$ smooth curve in the plane. Then the complex Fourier length spectrum $\{c_p, p\in\Z\}$ satisfies:
\begin{enumerate}
\item $c_0=L(\Gamma)>0$ and $\bar c_p=c_{-p}$. If $\Gamma$ is a closed (Jordan) curve then $c_{\pm1}=0$.
\item For any $N\in\N$, the Toeplitz circulant matrix $R=\hbox{Toep}(c_0,c_1, \ldots, c_{N-1})$, i.~e. $R_{pq}=c_{p-q}$, is a positive semidefinite complex Hermitian matrix.
\end{enumerate}
\end{tvr}

\prf
(i) We have 
\[
c_{\pm 1}=\int_\Gamma (t_1 \mp i t_2)  \ds = \int_\Gamma \partial_s (x_1 \mp i x_2)  \ds = 0
\]
because $\Gamma$ is a  closed curve and $(t_1,t_2)^T = \vect=\partial_s\vecx = (\partial_sx_1, \partial_s x_2)^T$. The rest of the statement (i) directly 
follows from the definition of the Fourier length spectrum.

In order to prove the statement (ii) we calculate
\begin{eqnarray*}
z^* R z &=&  \sum_{k,m=1}^N \bar z_k c_{k-m} z_m 
= \int_\Gamma \sum_{k,m=1}^N \bar z_k \exp(-i(k-m)\nu) z_m \,\ds 
\\
&=& \int_\Gamma \sum_{k,m=1}^N \bar z_k \exp(-i k \nu) z_m \exp(i m \nu)\,\ds
=\int_\Gamma \left| \sum_{k=1}^N  z_k \exp(i k \nu)\right|^2 \ds\ge 0,
\end{eqnarray*}
for any vector $z\in \CC^N$. Hence $R\succeq 0$, as claimed. \qed

\medskip
For a general positive semidefinite Toeplitz circulant complex matrix we can estimate off-diagonal terms by the diagonal ones. 

\begin{tvr}\label{Toeplitz-prop}
Let $c_0, c_1, \ldots, c_{N-1}\in \CC$. Assume that the complex Toeplitz circulant matrix  $R=\hbox{Toep}(c_0,c_1, \ldots, c_{N-1})$ is a positive semidefinite. Then
\begin{enumerate}
\item $|c_k|\le c_0$ for any $k\in \Z, |k|\le N-1$.

\item If, in addition, $c_1=0$ then $|c_{2k}|^2 + |c_{2k+1}|^2 \le c_0^2$ for any $k\le N/2-1$, and 
\begin{equation}
\sum_{p=2}^{N-1} \frac{|c_p|^2}{p^2-1} \le \frac{c_0^2}{2}\left(1-\frac{1}{N}\right).
\label{seriesestimate}
\end{equation}
\end{enumerate}

\end{tvr}

\prf Recall that an $N\times N$ Hermitian matrix $R$ is positive semidefinite if and only if  the main $K\times K$ submatrix 
$W$, $W_{pq} = R_{n_p n_q}$, is positive semidefinite for any index subset $\{n_1, \ldots, n_K\}\subseteq \{1,\ldots , N\}$ 
(see e.g. \cite[Theorem 6.2, p. 160]{Z}). To prove statement (i) it is sufficient to consider a $2\times2$ matrix $W$ 
corresponding to the index subset $\{1, k+1\}$, i.~e. 
$W=\left(\begin{array}{cc} c_0 & \bar c_k \\ c_k & c_0 \end{array}\right)$. Since $W\succeq 0$ we have $|c_k|\le c_0$ for each 
$k=1,\ldots, N-1$. 

In order to prove statement (ii) we consider the index subset $\{1, 2k+1, 2k+2\}$. Then the corresponding $3\times3$ matrix $W$ has the form 
\[W=\left(\begin{array}{ccc}
c_0 & \bar c_{2k} & \bar c_{2k+1}  \\ 
c_{2k} & c_0 & \bar c_1 \\ 
c_{2k+1} & c_1 & c_0
\end{array}\right)
=\left(\begin{array}{ccc}
c_0 & \bar c_{2k} & \bar c_{2k+1}  \\ 
c_{2k} & c_0 & 0\\ 
c_{2k+1} & 0 & c_0
\end{array}\right),
\]
because $c_1=0$. As $0\le\det(W)=c_0(c_0^2 - |c_{2k}|^2 - |c_{2k+1}|^2)$ we obtain the  estimate $|c_{2k}|^2 + |c_{2k+1}|^2 \le c_0^2$. 

To prove inequality (\ref{seriesestimate}) for $N$ even, we have
\begin{eqnarray*}
\sum_{p=2}^{N-1} \frac{|c_p|^2}{p^2-1}
&=& \sum_{k=1}^{N/2 -1} \frac{|c_{2k}|^2}{(2k)^2-1} + \frac{|c_{2k+1}|^2}{(2k+1)^2-1}
\\
&\le&
\sum_{k=1}^{N/2 -1} \frac{c_0^2 }{(2k)^2-1} = \frac{c_0^2}{2}\left(1-\frac{1}{N-1}\right).
\end{eqnarray*}

For $N$ odd, we can apply the above inequality as for an even dimension $N-1$ to obtain
\begin{eqnarray*}
\sum_{p=2}^{N-1} \frac{|c_p|^2}{p^2-1}
&=& 
\sum_{p=2}^{N-2} \frac{|c_p|^2}{p^2-1} +  \frac{|c_{N-1}|^2}{(N-1)^2-1}
\\
&\le& \frac{c_0^2}{2}\left(1- \frac{1}{N-2} +\frac{2}{(N-1)^2-1}\right)
=\frac{c_0^2}{2}\left(1-\frac{1}{N}\right).
\end{eqnarray*}
\qed

Applying Prop.~\ref{Toeplitz-prop} for the case of the Fourier length spectrum of a Jordan curve we obtain the following result.

\begin{dos}\label{Fspec-cor}
Let $\Gamma$ be a $C^1$ smooth Jordan curve in the plane. Then the complex Fourier length spectrum $\{c_p, p\in\Z\}$ satisfies inequality (\ref{seriesestimate}) for any $N\in \N$ and this estimate is optimal. 
\end{dos}

By optimality of the estimate we mean that there exists an $r$-parameterized family of Jordan curves for which the left hand side of (\ref{seriesestimate}) converges to $c_0^2/2$ as $N\to\infty$ and $r\to0$.

\prf
The proof follows directly from Prop.~\ref{Fspec-prop} and \ref{Toeplitz-prop}. To prove optimality of (\ref{seriesestimate}), let us consider a "capsule" like curve consisting of two horizontal line segments of the length $l>0$ connected by half-arcs with a radius $r>0$ from Example~\ref{fourier-spec-ex} (see Figure~\ref{fig:capsule}). Then $c_0=L(\Gamma) = 2l + 2\pi r$ and $c_{2k}=2 l, c_{2k+1}=0$. The left hand side of inequality (\ref{seriesestimate}) is therefore 
\[
\sum_{p=2}^{N-1} \frac{|c_p|^2}{p^2-1}=4 l^2 \sum_{k=1}^{N/2 -1} 1/((2k)^2-1)=2 l^2 (1-1/(N-1))
\]
and it tends to the value $c_0^2/2$ as $N\to\infty$ and $r\to0$. 
\qed

In subsequent sections we will prove that inequality (\ref{seriesestimate}) plays an essential  role in the proof of the fact that the enhanced semidefinite relaxation method for solving the inverse Wulff problem indeed yields the optimal solution for the anisoperimetric function $\sigma$.

\begin{tvr}\label{pozdef}
Let $P_0=\hbox{diag}(p_0, p_1, \ldots, p_{N-1}, q_0, q_1, \ldots, q_{N-1})$ be an $n\times n$ real matrix, $n=2N$. Assume $p_0<0,q_0\le 0, p_1=q_1=0$, and $p_k,q_k>0$, for $k=2,\ldots, N-1$. Let $A$ be a $2\times n$ real matrix $A=\left(
\begin{array}{cc}
 \alpha^T & \beta^T \\ 
 0_N^T & e_1^T
\end{array} 
\right)
$
where $\alpha=(\alpha_0, \alpha_1, \alpha_2, \ldots, \alpha_{N-1})^T$, $\beta=(\beta_0, \beta_1, \beta_2,\ldots, \beta_{N-1})^T \in\R^N$ are such that $\alpha_1=\beta_1=0$. Assume $\varrho >  -q_0\ge 0$ and 
\begin{equation}
\sum_{k=2}^{N-1}\left(\frac{\alpha_k^2}{p_k}+\frac{\beta_k^2}{q_k}\right)
\leq \frac{\alpha_0^2}{-p_0} - \frac{1}{\varrho}-\frac{\beta_0^2}{q_0+\varrho}.
\label{sum-cond}
\end{equation}
Then the matrix $P_0 + \varrho A^T A$ is positive semidefinite.

\end{tvr}

\prf
We note that that we can delete zero columns and rows from the matrix $P_0+\varrho A^TA$. Since $\alpha_1=\beta_1=0$ the matrix $P_0 + \varrho A^T A \succeq 0$ if and only if the squeezed $(n-2)\times(n-2)$ matrix $\tilde P_0+\varrho \tilde A^T\tilde A\succeq 0$ is positive semidefinite in which the second and $N+2$ zero columns and rows of $P_0 + \varrho A^T A$ were omitted. The matrix $\tilde P_0+\varrho \tilde A^T\tilde A$ has the following structure:
\[
\tilde P_0+\varrho \tilde A^T\tilde A = \left(
\begin{array}{cc}
 p_0 +\varrho \alpha_0^2  & \varrho \alpha_0 v^T \\ 
 \varrho \alpha_0 v  & D +\varrho v v^T
\end{array} 
\right),
\]
where $v= (\alpha_2,\ldots, \alpha_{N-1}, \beta_0, \beta_2, \ldots, \beta_{N-1})^T\in\R^{N-3}$. The main diagonal submatrix $D=\hbox{diag}(p_2, \ldots, p_{N-1}, q_0+\varrho, q_2, \ldots, q_{N-1})$ is positive definite provided that $\varrho >- q_0\ge 0$. Hence the block submatrix $D+\varrho v v^T\succ 0$. By using the Schur complement property, we conclude that $\tilde P_0+\varrho \tilde A^T\tilde A  \succeq 0$ if and only if the Schur complement is nonnegative, i.~e. 
\begin{equation}
0\le p_0+\varrho \alpha_0^2 - \varrho^2 \alpha_0^2\,  v^T(D +\varrho v v^T )^{-1} v.
\label{shurcond}
\end{equation}
By using the Morrison-Sherman formula we obtain $(D +\varrho v v^T )^{-1} = D^{-1} -\frac{\varrho}{1+\varrho\gamma} D^{-1} v v^T D^{-1}$, where we have denoted $\gamma=v^T   D^{-1} v \ge 0$. Hence condition (\ref{shurcond}) is equivalent to the inequality:
\[
0\leq p_0+\varrho \alpha_0^2 - \varrho^2 \alpha_0^2 \left( \gamma - \frac{\varrho \gamma^2}{1+\varrho\gamma}\right)
= \frac{p_0 + \varrho p_0 \gamma + \varrho \alpha_0^2}{1+\varrho\gamma}.
\]
Since $p_0<0$ then solving the above inequality for $\gamma=v^T   D^{-1} v$ yields the condition $v^T   D^{-1} v \le \alpha_0^2/(-p_0) - 1/\varrho$, which is indeed condition (\ref{sum-cond}).

\qed

\section{Enhanced semidefinite relaxation method}

\subsection{General nonconvex quadratic optimization problem with linear matrix inequality constraints}

In this section, our goal is to propose and investigate the so-called enhanced semidefinite relaxation method  for solving the following optimization problem:
\begin{equation}\label{P}
\begin{array}{rl}
\min & x^TP_0x+2q_0^Tx+r_0 \\
{\rm s. t.} & x^T P_l x+2q_l^Tx+r_l \leq 0, \quad l=1,\ldots, d, \\  
 & Ax=b, \\
 & H_0+\sum_{j=1}^n x_j H_j \succeq 0,
\end{array} 
\end{equation}
where $x\in\mb{R}^n$ is the variable and the data: $P_0, P_l \in \mathcal{S}^n$ are $n\times n$ real symmetric matrices, $q_0, q_l\in\mb{R}^n$, $r_0, r_l\in\mb{R}$, $A$ is an $m\times n$ real matrix, $b\in\mb{R}^m$ and $H_0, H_1, \ldots, H_n\in \mathcal{H}^k$, i.~e. $H_0, H_1, \ldots, H_n$ are $k\times k$ complex Hermitian matrices. The last constraint in \eqref{P} is a complex linear matrix inequality (LMI). It can be easily transformed into  real LMI using the following equivalence
\begin{equation}\label{realvscomplex}
H\succeq 0 \ \Longleftrightarrow \ 
\tilde H \equiv \left(\begin{array}{lr}
\Re H & - \Im H \\
\Im H & \Re H
\end{array}\right)\succeq 0.
\end{equation}
Regarding the input matrices $P_0, P_l, l=1,\ldots, d,$ and $A$ we will henceforth assume the following assumption:

\begin{itemize}
\item[\textbf{(A)}] $P_l\succeq 0$ for $l=1,\ldots, d$ and there exists a real $m\times n$ matrix $V$ such that
\[
P_0 + \frac12 (V^T A+A^T V)\succeq 0.
\]
\end{itemize}

\begin{rem}\label{rem:finsler}
Assumption (A), in particular the condition $P_0+\frac12 (V^T A+A^T V)\succeq 0$ for some  $V$ includes a special case when $P_0$ is positive semidefinite on the null space $\{x \ | \ Ax=0\}$. In such a case, it follows from the Finsler theorem (c.f. \cite{finsler}) that the matrix $P_0+\varrho A^T A\succeq 0$ for each $\varrho>0$ sufficiently large. Assumption (A) is then satisfied if we set $V=\varrho A$. However, assumption (A) is more general. Indeed, let us consider the following example: 
$P_0=\left(\begin{array}{rr}
0 & -1 \\
-1 & 1 \end{array}\right), \quad A= \left(\begin{array}{cc} 0 & 1\\
\end{array}\right)$.
Then there exists no real number $\varrho\ge 0$ such that $P_0+\varrho A^T A = \left(\begin{array}{rr}
0 & -1 \\
-1 & 1+\varrho
\end{array}\right)\succeq 0$.
However, for the choice of $V= \left(\begin{array}{cc} 2 & 0\\
\end{array}\right)$ we have
$
P_0+\frac12 (V^T A+A^T V)=
\left(\begin{array}{rr}
0 & 0 \\
0 & 1
\end{array}\right)\succeq 0.$
\end{rem}

In what follows,  under assumption (A), we will show that problem \eqref{P} can be solved by means of Lagrangian duality and relaxation with a convex semidefinite programming (SDP) problem. Note that for the case $A=0, H_j=0, j=0,1,\ldots, n$, and $d=1$ the method was proposed and analyzed in \cite[Appendix C.3]{bova}. In what follows, we will propose a relaxed SDP that includes LMI of the general form $H_0+\sum_{j=1}^n x_j H_j \succeq 0$. Moreover, we will augment SDP  \eqref{P} with the additional quadratic-linear constraint $Ax x^T=bx^T$ which follows from $Ax=b$. 

\subsection{Augmented problem and enhanced semidefinite relaxation}

Clearly, problem \eqref{P} is equivalent to the following augmented problem with one additional constraint:
\begin{equation}\label{P1}
\begin{array}{rl}
\min & x^TP_0x+2q_0^Tx+r_0 \\
\hbox{s. t.} & x^T P_l x+2q_l^Tx+r_l \leq 0, \quad l=1,\ldots, d, \\  
 & Ax=b, \quad  Ax x^T=bx^T,\\
 & H_0+\sum_{j=1}^n x_j H_j \succeq 0. 
\end{array} 
\end{equation}
We will show that the additional constraint $Ax x^T=bx^T$ becomes a linear constraint $A X = b x^T$ between the relaxed matrix $X\succeq x x^T$ and the vector $x$. Furthermore, we will prove that the original problem and its second Lagrangian dual yield the same optimal values provided that the matrices $P_l, l\ge 0$, and $A$ satisfy assumption (A). In particular, if $P_0+\varrho A^T A\succeq 0$ for $\varrho \gg1$,  then, with regard to Remark~\ref{rem:finsler} the value function  $x\mapsto x^TP_0x+2q_0^Tx+r_0$ is convex on the affine subspace $\{x\ |\ Ax=b\}$ of the feasible set of the  semidefinite relaxed problem. This is why we will henceforth refer a method when additional constraint $Ax x^T=bx^T$ is added to $Ax=b$ to as {\it the enhanced semidefinite relaxation method}. 

The idea of semidefinite relaxation of \eqref{P1} is rather simple and it consists in relaxing the equality $X=x x^T$ by the semidefinite inequality $X\succeq x x^T$. Although the form of the relaxed problem can be deduced from \eqref{P1} we will still  present a systematic way of its derivation based on construction of the second Lagrangian dual SDP to \eqref{P1}. 

\subsection{The first and second Lagrangian dual problems}

Next, we construct the first and second Lagrange dual problem for the augmented problem \eqref{P1}. To this end, let us consider the following Lagrangian function ${\mathcal L}^1={\mathcal L}^1(x; \lambda, u, V, Z)$:
\begin{eqnarray*}
{\mathcal L}^1 &=& x^T P_0x+2q_0^Tx+r_0 + \sum_{l=1}^d \lambda_l [x^T P_l x+2q_l^Tx+r_l] + u^T(Ax-b) \\
&& + \hbox{tr} (V^T(A x x^T-b x^T)) - \hbox{tr} (Z^T(\tilde{H}_0+\sum_{j=1}^n x_j \tilde{H}_j)),
\end{eqnarray*}
where $0\leq \lambda\in\R^d, u\in \mb{R}^n$, $V$ is an $m\times n$ real matrix and $Z$ is a $2k\times 2k$ real symmetric positive semidefinite matrix. The tuple $(\lambda,u,V,Z)$ represents the Lagrange multipiers to problem \eqref{P1}. Here we have used a real version of complex LMI based on equivalence \eqref{realvscomplex}. The dual problem can be obtained by analyzing $\inf_x{\mathcal L}^1(x; \lambda, u, V, Z)$. In Appendix we show by using straightforward calculations and applying properties of the Schur complement, the Lagrangian dual problem \eqref{P1} has the form: 
\begin{equation}\label{D2}
\begin{array}{rl}
\max & \gamma\\
\hbox{s. t. }
  & M_0+\sum_{l=1}^d \lambda_l M_l + M_*(V,u)-\sum_{j=0}^nz_j N_j -\gamma N_0\succeq 0,\\
  & Z\succeq 0, \ \lambda\geq 0, \\
  & z_j=\hbox{tr}(Z^T\tilde{H}_j),\ j=0,1,\ldots, n.
\end{array}
\end{equation}
Here we have denoted 
\[
M_j=\left(\begin{array}{cc}
P_l & q_l \\
q_l^T & r_l
\end{array}\right),\ M_*(V,u)=\frac12
\left(\begin{array}{cc}
V^TA+A^TV & A^Tu-V^Tb \\
u^TA-b^TV & -2u^Tb
\end{array}\right)
,
\]
\[
N_0=\left(\begin{array}{cc}
0_{n\times n} & 0_n \\
0_n^T & 1
\end{array}\right), \quad
N_j=\frac12 \left(\begin{array}{cc}
0_{n\times n} & e_j \\
e_j^T & 0
\end{array}\right),
\]
where $e_j$ is the $j$-th unit vector in $\mb{R}^n$. Note that problem \eqref{D2} is closely related to the so-called Shor-relaxation method (see \cite{shor} for details).

\medskip
We proceed by constructing the second Lagrangian dual problem. Let us consider problem \eqref{D2}. We define its Lagrangian ${\mathcal L}^2(\gamma,\lambda, Z,V,u,z; W,\beta,\tilde{X},\alpha)$ function as follows:
\begin{eqnarray*}
{\mathcal L}^2 &=&\gamma+\hbox{tr}(ZW)+\lambda\beta+\hbox{tr}(\tilde{X}(M_0+ \sum_{l=1}^d \lambda_l M_l + M_*(u,V) 
\\
&&-\sum_{j=0}^n z_jN_j-\gamma N_0))+\sum_{j=0}^n\alpha_j(z_j-\hbox{tr}(Z\tilde{H}_j))
\end{eqnarray*}
with dual variables $W\succeq 0, \tilde{X}=\left(\begin{array}{cc}
X & x \\ x^T & \varphi
\end{array}\right)\succeq 0, \ \beta\geq 0, \ \alpha_j\in\mb{R}, \ j=0,1,\ldots, n$. 
By $\beta\ge0$ we mean $\beta_j\ge 0$ for each $j=0,\dots, N-1$. The dual problem to \eqref{D2} can be obtained by solving the following problem
\[
\sup_{\gamma,\lambda, Z,V,u,z} {\mathcal L}^2(\gamma,\lambda, Z,V,u,z; W,\beta,\tilde{X},\alpha).
\]
After straightforward calculations (see Appendix for details) the second Lagrangian dual to SDP \eqref{P1} then reads as follows:
\begin{equation}\label{D3}
\begin{array}{rl}
\min & \hbox{tr}(P_0 X)+2q_0^Tx+r_0\\
\hbox{s. t. } & \hbox{tr}(P_l X)+2q_l^Tx+r_l\leq 0, \ \ l=1,\ldots,d,\\
    & Ax=b, \ AX=bx^T, \  X\succeq x x^T,\\
& {H}_0+\sum_{j=1}^n x_j {H}_j\succeq 0.
\end{array}
\end{equation}

Henceforth, we will refer \eqref{D3} to as the enhanced semidefinite relaxation of problem \eqref{P}. Notice that the second Lagrangian dual to \eqref{P} is just problem \eqref{D3} without the constraint $AX=b x^T$.

\begin{rem}\label{rem-DD3}
Instead of the additional constraint $Axx^T=bx^T$ in \eqref{P1} we could alternatively use the simplified constraint $\hbox{tr}(A^TAxx^T)=x^TA^TAx=b^TAx$ and consider the augmented problem \eqref{P1} in which the constraint $A x x^T = b x^T$ is replaced by the equation $x^TA^TAx=b^TAx$. Using a similar technique as before we can construct the Lagrange dual which is just optimization problem \eqref{D2} with $V=\nu A$ where $\nu\in \R$. 
Then the second Lagrangian dual problem has the form of optimization problem \eqref{D3} in which the constraint $AX=bx^T$ is replaced by $\hbox{tr}(A^TAX)=b^TAx$. In such a case,  assumption (A) has to be modified by taking constraint $V=\varrho A$. 
\end{rem}

Problem \eqref{D3}  can be viewed as the semidefinite relaxation \eqref{P1}. Replacing the condition $X\succeq xx^T$ in \eqref{D3} with $X=xx^T$ would lead to a problem equivalent to problem \eqref{P1}. Such a relaxation is often used in solving nonconvex quadratic problems or combinatorial optimization problems (see \cite{bovasdp}, \cite{review}, \cite{nowak}, \cite{sdprelax}).

\begin{rem}
It is worth noting that problem \eqref{D3} has no interior point unless $A=0$. Indeed, suppose that there are $X$ and $x$ feasible for \eqref{D3} and satisfying  $X\succ xx^T$. So there exists a positive definite matrix $D$ such that $X=xx^T+D$. But then $bx^T=AX=Axx^T+AD=bx^T+AD$. Hence $AD=0$ and since $D$ is nonsingular, we obtain $A=0$. Similar property holds for problem \eqref{D3} in which the constraint $AX=bx^T$ is replaced by $\hbox{tr}(A^TAX)=b^TAx$ (see Remark~\ref{rem-DD3}). Indeed, $b^TAx=\hbox{tr}(A^TAX)=x^TA^TAx+\hbox{tr}(A^TAD)$ implies $\hbox{tr}(A^TAD)=0$. However, since $D\succ 0$, we obtain $A=0$.
\end{rem}
 
\subsection{Equivalence of problems}

In this section we give sufficient conditions under which problem \eqref{P1} and its second dual \eqref{D3} yield the same optimal values. 

\begin{tvr}\label{tvr1}
Denote $\hat p_1$ the optimal value of problem \eqref{P} and $\hat p_2$ the optimal value of enhanced semidefinite relaxed problem  \eqref{D3}. Then $\hat p_1\geq \hat p_2$.
\end{tvr}

\prf
If problem \eqref{P1} is not feasible, then $\hat p_1=+\infty$ and the inequality is satisfied.
Assume $\hat p_1<+\infty$. Denote $f(x)=x^TP_0x+2q_0^Tx+r_0$ the quadratic objective of  problem \eqref{P1}. Then there exists a sequence of points $\{x^{(k)}\}_{k=1}^{\infty}$, feasible for \eqref{P1} such that $\hat p_1=\liminf_{k\to\infty} f(x^{(k)})$. It covers the case when optimum is attained as well as optimum is not attained. Denote $X^{(k)}=x^{(k)} (x^{(k)})^T\succeq 0, k=1,2,\ldots, \infty$. Then it is easy to see that the pair $(x^{(k)}, X^{(k)})$ is feasible for the problem \eqref{D3} for all $k=1,2,\ldots, \infty$. Since $x^TP_j x=\hbox{tr}(P_j x x^T), j=0,1, \ldots, d$, we have $f(x^{(k)})\geq \hat p_2$ for each $k\ge 1$. Hence  $\hat p_1\geq \hat p_2.$ \qed

If we denote $\hat d$ the optimal value of problem \eqref{D2}, the previous proposition together with the weak duality yields $\hat d\leq \hat p_2\leq \hat p_1$. A strong duality property between problems \eqref{P1} and \eqref{D2} would imply the equality $\hat p_1=\hat p_2$. 
In the case $d=1, A=0, b=0$ and there are no LMI constraints (i.~e. $H_i=0$) in \eqref{P} then  the strong duality has been shown under the assumption that problem \eqref{P} has an interior point (c.f. \cite[Appendix B.1]{bova}). The proof relies on the so-called S-procedure and it is not obvious how to generalize it to the case of nontrivial LMI constraints and/or quadratic-linear constraints occurring in \eqref{P1}. Nevertheless, in the following proposition we prove the equality $\hat p_1=\hat p_2$ under assumption (A) made on input matrices $P_1, P_l$ without assuming the strong duality property. First we introduce an auxiliary lemma whose proof easily follows from the property of positive semidefinite matrices: $\hbox{tr}(AB)\geq 0$ for $A\succeq 0, B\succeq 0$.

\begin{lem}\label{lema1}
Let $M\succeq 0$ and $X\succeq xx^T$. Then $\hbox{tr}(MX)\geq x^TMx$.
\end{lem}

\begin{vet}\label{prop-ekviv}
Suppose that the SDP problem \eqref{P} is feasible and assumption (A) is satisfied. Let $\hat p_1$ be the optimal value of \eqref{P} and $\hat p_2$ be the optimal value of SDP \eqref{D3} obtained by the enhanced semidefinite relaxation method. Then $\hat p_1=\hat p_2$. If $(\tilde x, \tilde X)$ is an optimal solution to  \eqref{D3} then $\tilde x$ is the optimal solution to \eqref{P}.
\end{vet}

\prf
Let $(x,X)$ be a feasible solution \eqref{D3}. By means of Lemma \ref{lema1} we have
\[
x^T P_l x + 2 q_l^T x + r_l \le  \hbox{tr}(P_l X) + 2 q_l^T x + r_l\le 0, \quad l=1,\ldots, d.
\] 
Hence $x$ is a feasible solution to \eqref{P1}. Since $P_0+ \frac12 (V^T A + A^T V) \succeq 0$ and 
\[
\hbox{tr}((V^T A + A^T V) (X - x x^T)) = 2 \hbox{tr}(V^T ( AX - Ax x^T) ) = 0
\]
for any $(x,X)$ feasible to \eqref{D3}. By Lemma~\ref{lema1} we furthermore have $x^T P_0 x  + 2 q_0^T x + r_0 \le \hbox{tr}(P_0 X) + 2 q_0^T x + r_0 \equiv \Phi(x,X)$. 

In order to prove the equality $\hat p_1=\hat p_2$ we consider a minimizing sequence $(x^{(k)}, X^{(k)})$ of feasible solutions to \eqref{D3},  $\hat p_2=\liminf_{k\to\infty} \Phi(x^{(k)}, X^{(k)})$.  
Since $x^{(k)}$ is feasible to \eqref{P1} we have $\hat p_1 \le (x^{(k)})^T P_0 x^{(k)}  + 2 q_0^T x^{(k)} + r_0 \le  \Phi(x^{(k)}, X^{(k)})$ 
for any $k\in\N$.  Thus $\hat p_1 \le \liminf_{k\to\infty} \Phi(x^{(k)}, X^{(k)}) = \hat p_2$. Finally, if  $(\tilde x, \tilde X)$ 
is an optimal solution to  \eqref{D3} then 
$\hat p_1 \le \tilde x^T P_0 \tilde x  + 2 q_0^T \tilde x + r_0 \le \hbox{tr}(P_0 \tilde X) + 2 q_0^T \tilde x + r_0 =\hat p_2$. 
Hence $\tilde x$ is optimal to \eqref{P1}, as claimed.

\qed

\subsection{Application of the Enhanced Semidefinite Relaxation Method to a solution of the inverse Wulff problem}

We conclude this section with construction of the second Lagrangian dual formulation of optimization problem \eqref{maxWNP} resolving minimal anisoperimetric ratio over all anisotropy function belonging to the cone $\K^N$. 

\begin{vet}
The enhanced semidefinite relaxation of optimization problem \eqref{maxWNP} has the form
\begin{equation}\label{maxWNP-relax}
\begin{array}{rl}
\min & \hbox{tr}( P_0 X) \\
{\rm s. t.} & A x = b,\  A X = b x^T, \ X \succeq x x^T,  \ x=[x^R; x^I],
\\ & 
\sum_{p=k+1}^N F_{p,p-k} = x^R_k + i x^I_k, 
\\ &
\sum_{p=k+1}^N G_{p,p-k} = (1-k^2)(x^R_k + i x^I_k),
\\
& \hbox{for}\ k=0,\ldots, N-1,
\\ 
& F, G\succeq 0,
\end{array} 
\end{equation}
where the matrices $P_0, A$ and $b$ are defined as in \eqref{P0def} and \eqref{Abdef}. Problem \eqref{maxWNP} is feasible. Optimal values $\hat p_1, \hat p_2$ of \eqref{maxWNP} and \eqref{maxWNP-relax}, respectively, are finite and $\hat p_1 = \hat p_2$. If  $(\tilde x,\tilde X)$ is an optimal solution to \eqref{maxWNP-relax} then $\tilde x$ is the optimal solution to \eqref{maxWNP}. Conversely, if $\tilde x$ is the optimal solution to \eqref{maxWNP} then $(\tilde x, \tilde X)$ is the optimal solution to \eqref{maxWNP-relax} where $\tilde X=\tilde x \tilde x^T$.
\end{vet}

\prf Feasibility of \eqref{maxWNP} is obvious because for $\tilde\sigma=(1,\ldots,0)^T\in\K_N$ we have  $L_{\tilde\sigma}(\Gamma)=L(\Gamma)$ and $|W_{\tilde\sigma}|=\pi>0$. Furthermore, we have the following estimate:
\[
- x^T P_0 x = |W_\sigma| = \frac{L_\sigma(\Gamma)^2}{4\Pi_\sigma(\Gamma){\mathcal A}(\Gamma)}\le \frac{L(\Gamma)^2}{4 {\mathcal A}(\Gamma)} <\infty, 
\]
for any $\sigma = x^R + i x^I \in\K^N,\ x=[x^R; x^I],$ feasible to \eqref{maxWNP}. So the optimal value $\hat p_1$ of SDP \eqref{maxWNP} is finite. 

Note that problem \eqref{maxWNP} can be rewritten in the form of SDP \eqref{P}. Since the quadratic constraints become linear when assuming $P_l \equiv0$ for $l=1,\ldots,d$ we can construct the second Lagrangian dual to the augmented SDP \eqref{P1} having additional linear constraints of the form $Qx + r =0$. Furthermore, by taking a standard basis of ${\mathcal H}^N$, the semidefinite constraints $F,G\succeq 0$ can rewritten as LMI in the form: ${H}_0+\sum_{j=1}^n x_j {H}_j\succeq 0$. 

Taking into account Prop.~\ref{Toeplitz-prop}, for $\alpha_k+i \beta_k = 2 c_k, k\ge 1$, $\alpha_0=c_0, \beta_0=0$ and $p_0=q_0=-\pi, p_k=q_k= 2\pi (k^2-1)$ we conclude that condition (\ref{sum-cond}) is fulfilled for all $\varrho> \pi \max(N/c_0^2,1)$. Hence $P_0 + \varrho A^T A \succeq 0$ and assumption (A) is satisfied. By Theorem~\ref{prop-ekviv} we conclude $\hat p_1 = \hat p_2$ where $\hat p_2$ is the optimal value of the SDP \eqref{maxWNP-relax} obtained by the enhanced semidefinite relaxation method. 
\qed

\section{Convergence analysis}
In this section we prove convergence of a sequence of approximative anisoperimetric ratio to the optimal value of problem \eqref{anisoproblem}. For any finite dimension $N\in\N$ we recall that $\sigma^N \in\K^N$ is a minimizer of the $N$-dimensional restriction \eqref{maxWN} of the original problem \eqref{anisoproblem}. Then, for  $\tilde\sigma^{N+1} =  (\sigma^N_0, \sigma^N_1, \ldots, \sigma^N_{N-1}, 0)^T\in\K^{N+1}$ we have $L_{\tilde\sigma^{N+1}}(\Gamma)= L(\Gamma)$ and this is why $\tilde\sigma^{N+1}$ is feasible solution to \eqref{maxWN} in the dimension $N+1$. Thus we obtain
$|W_{\sigma^N}| = |W_{\tilde\sigma^{N+1}}| \le  |W_{\sigma^{N+1}}| \quad\hbox{for all}\ N\in\N$. It means that $1\le \Pi_{\sigma^{N+1}}(\Gamma) \le \Pi_{\sigma^{N}}(\Gamma)$ for each $N\in\N$. This is why the sequence $\{ \Pi_{\sigma^{N}}(\Gamma)\}_{N=1}^\infty$ of anisoperimetric ratios is nonincreasing and having thus a finite limit. More precisely, we have the following result:

\begin{vet}\label{convergence}
Let $\sigma^N \in\K^N$ be a minimizer to optimization problem \eqref{maxWN} in the dimension $N\in\N$. Then
$1\le \lim_{N\to\infty} \Pi_{\sigma^{N}}(\Gamma) = \inf_{\sigma\in\K} \Pi_\sigma(\Gamma)$.
\end{vet}

\prf
Let $\sigma\in\K$ be fixed and such that $\Pi_\sigma(\Gamma)<\infty$. Given the dimension $N\in\N$ we will construct $\tilde\sigma^N\in\K^N$ such that $\tilde\sigma^N \to\sigma$ as $N\to\infty$ in the norm of the Sobolev space $W^{1,2}_{per}(0,2\pi)$. To this end, we employ both the Bochner and McLean criteria for positiveness of Fourier series (c.f. Prop.~\ref{Dumit-pos} and \ref{Lean-pos}). 

Since $\sigma\in\K$ it follows from Prop.~\ref{Dumit-pos} that the Toeplitz matrices
\[
Q^{(N)} = \hbox{Toep}(\sigma_0, \sigma_1, \ldots, \sigma_{N-1}), \ \ S^{(N)} = \hbox{Toep}(\xi_0, \xi_1, \ldots, \xi_{N-1}), 
\]
where $\xi_k:=(1-k^2)\sigma_k$, are positive semidefinite. Hence the sets $\mathcal{F}, \mathcal{G}\subset \mathcal{H}^N$ given by:
\begin{eqnarray*}
\mathcal{F}&=&\{ F \ |\  F\succeq 0,\ \sum_{p=k+1}^N F_{p,p-k} = \frac{N-k}{N} \sigma_k, \ \forall\ k=0, \ldots, N-1\},
\\
\mathcal{G}&=&\{ G\ | \  G\succeq 0,\ \sum_{p=k+1}^N G_{p,p-k} = \frac{N-k}{N} (1-k^2)\sigma_k, \ \forall\  k=0, \ldots, N-1\}
\end{eqnarray*}
are nonempty as $\frac1N Q^{(N)}\in\mathcal{F}$ and $\frac1N S^{(N)}\in\mathcal{G}$, respectively. By Prop.~\ref{Lean-pos} we have 
\[ 
\tilde\sigma^N \in \K^N, \quad \hbox{where} \ \ \tilde\sigma^N_k = \frac{N-k}{N} \sigma_k, \ \  k=0, \ldots, N-1.
\]
Then the distance between $\tilde\sigma^N$ and $\sigma$ in the norm of the Sobolev space $W^{1,2}_{per}(0,2\pi)$ can be estimated as follows:
\begin{eqnarray*}
\Vert \tilde\sigma^N - \sigma\Vert_{1,2}^2 
&=&  \sum_{k=0}^{N-1} (1+k^2) |\tilde\sigma^N_k - \sigma_k|^2
+ \sum_{k=N}^{\infty} (1+k^2) |\sigma_k|^2
\\
&=&  \sum_{k=0}^{N-1} (1+k^2) \frac{k^2}{N^2} |\sigma_k|^2
+ \sum_{k=N}^{\infty} (1+k^2) |\sigma_k|^2  
\\
&\le& \frac{1}{N^2}\sum_{k=0}^{\infty} (1+k^2) k^2  |\sigma_k|^2 \le \frac{2}{N^2}\sum_{k=0}^{\infty} (1+k^4) |\sigma_k|^2 = \frac{2}{N^2}\Vert\sigma\Vert_{2,2}^2.
\end{eqnarray*}
As $\sigma\in\K\subset W^{2,2}_{per}(0,2\pi)$ we have $\Vert\sigma\Vert_{2,2}<\infty$ and therefore $\lim_{N\to\infty} \tilde\sigma^N = \sigma$ in the norm of $W^{1,2}_{per}(0,2\pi)$. More precisely, $\Vert \tilde\sigma^N - \sigma\Vert_{1,2}= O(N^{-1})$ as $N\to\infty$.

Clearly, the Wulff shape area $|W_\sigma|$ as well as the total interface energy $L_\sigma(\Gamma)$ are continuous functionals in $\sigma$ in the norm of the Sobolev space  $W^{1,2}_{per}(0,2\pi)$. Hence 
$\lim_{N\to\infty} |W_{\tilde\sigma^N}| = |W_\sigma|$ and $\lim_{N\to\infty} L_{\tilde\sigma^N}(\Gamma) = L_\sigma(\Gamma)$. Thus 
\[
\lim_{N\to\infty} \Pi_{\tilde\sigma^N}(\Gamma) = \Pi_\sigma(\Gamma).
\]
Finally, as $\tilde\sigma^N\in\K^N$ and $\sigma^N$ is a minimizer of $\Pi_\sigma$ in $\K^N$ we have $\Pi_{\tilde\sigma^N}(\Gamma)\ge \Pi_{\sigma^N}(\Gamma)$. Therefore $1\le \lim_{N\to\infty} \Pi_{\sigma^N}(\Gamma) \le \Pi_\sigma(\Gamma)$ for any $\sigma\in\K$ and the proof follows.
\qed

\begin{rem}\label{rem-capsule}
In the statement of Theorem~\ref{convergence} the infimum $\inf_{\sigma\in\K} \Pi_\sigma(\Gamma)$ need not be attained by any $\sigma\in\K$. Indeed, let us consider a convex "capsule" curve $\Gamma$ from Example~\ref{fourier-spec-ex-capsule}. Then $\inf_{\sigma\in\K} \Pi_\sigma(\Gamma)=1$ because the boundary $\partial W_\sigma$ of the limiting optimal Wulff shape should coincide with the convex curve $\Gamma$. If there is $\sigma\in\K\subset W^{2,2}_{per}(0,2\pi)$ such that $\Pi_\sigma(\Gamma)=1$ then for the curvature of $\partial W_\sigma=\Gamma$ we have $\kappa=[\sigma+\sigma'']^{-1} \ge 0$. Since  $\partial_s\nu = \kappa$, i.~e. $\dnu=\kappa\ds$, we obtain 
\begin{equation}
\int_{\partial W_\sigma} \frac{1}{\kappa} \ds = \int_0^{2\pi} \frac{1}{\kappa^2} \dnu = \int_0^{2\pi} [\sigma(\nu)+\sigma''(\nu)]^2 \dnu < \infty, 
\label{recip}
\end{equation}
because $\sigma$ and $\sigma''$ are square integrable functions for any $\sigma \in W^{2,2}_{per}(0,2\pi)$. But the curvature $\kappa\equiv0$ on nontrivial line segments of the capsule $\partial W_\sigma=\Gamma$. So $\int_{\partial W_\sigma} \frac{1}{\kappa} \ds = \infty$, a contradiction to \eqref{recip}.

\end{rem}

\section{Numerical experiments}

Let $\Gamma$ be a Jordan curve in the plane $\R^2$ and $\vecx^{(0)}, \vecx^{(1)}, \ldots, \vecx^{(K)}\in\Gamma$ be a set of its points where $\vecx^{(0)}=\vecx^{(K)}$. The curve $\Gamma$ will be approximated by a polygonal curve with vertices $\vecx^{(0)}, \vecx^{(1)}, \ldots \vecx^{(K)}$. The unit tangent $\vect^{(k)}$ vector at $\vecx^{(k)}$ will be approximated by 
$\vect^{(k)} \equiv (\vecx^{(k+1)}-\vecx^{(k-1)})/\Vert \vecx^{(k+1)}-\vecx^{(k-1)}\Vert$. Since $\ds =\Vert\partial_u \vecx\Vert \du \approx \frac12 \Vert \vecx^{(k+1)}-\vecx^{(k-1)}\Vert$ the elements of the Fourier length spectrum $\{c_p, p\in \Z\}$ can be approximated by
\begin{equation}
c_p = \int_\Gamma (t_1 - i  t_2\bigl)^p\ds \approx \frac12 \sum_{k=1}^{K-1} \bigr(t^{(k)}_1 - i  t^{(k)}_2\bigl)^p \Vert \vecx^{(k+1)}-\vecx^{(k-1)}\Vert.
\label{cp-discr}
\end{equation}

The enhanced semidefinite relaxation \eqref{maxWNP-relax} of optimization problem \eqref{maxWNP} was solved by using the powerful nonlinear convex programming solver SeDuMi developed by J.~Sturm \cite{sturm}. SeDuMi (Self-Dual-Minimization) implements self-dual embedding method proposed by Ye, Todd and Mizuno \cite{todd}. It is implemented as  an add-on for MATLAB and it has a capacity in solving large optimization problems, including \eqref{maxWNP-relax}. Without assuming the quadratic-linear constraint $AX=b x^T$ with $X\succeq x x^T$ in \eqref{maxWNP-relax} the SeDuMi solver was unable to solve the problem because of its unboundedness.

\begin{figure}
\begin{center}
\subfloat[]{\includegraphics[width=0.28\textwidth]{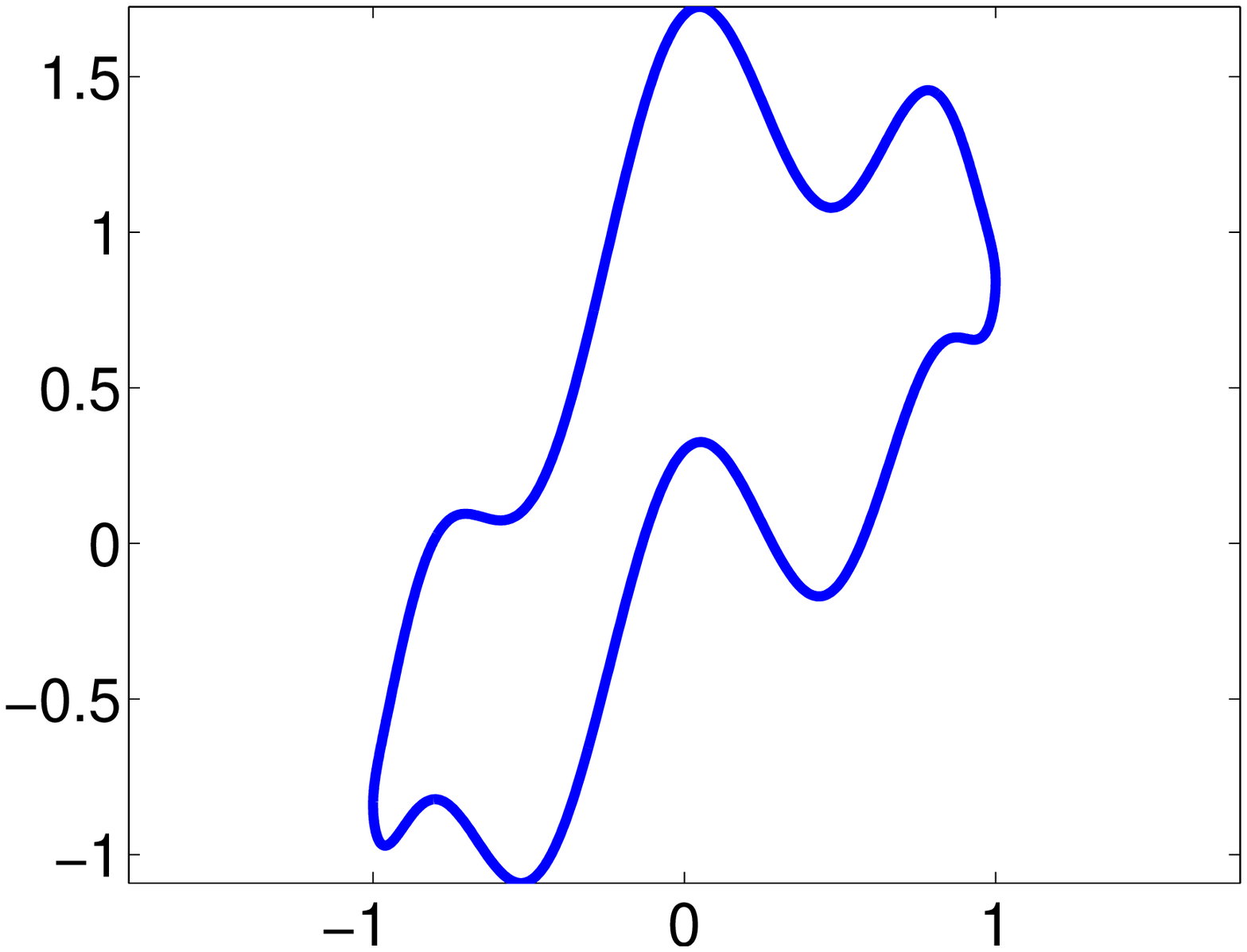}}
\subfloat[]{\includegraphics[width=0.28\textwidth]{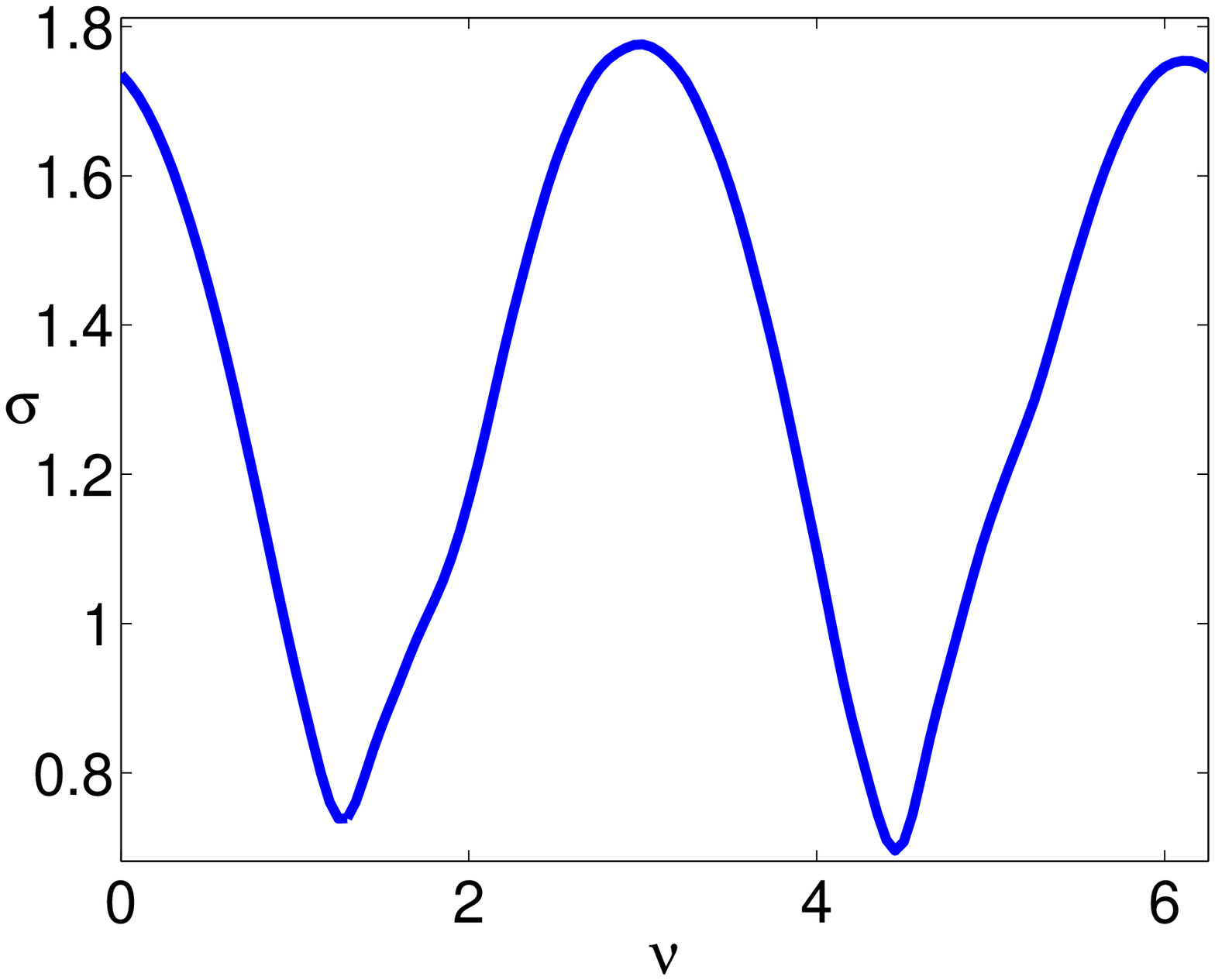}}
\subfloat[]{\includegraphics[width=0.28\textwidth]{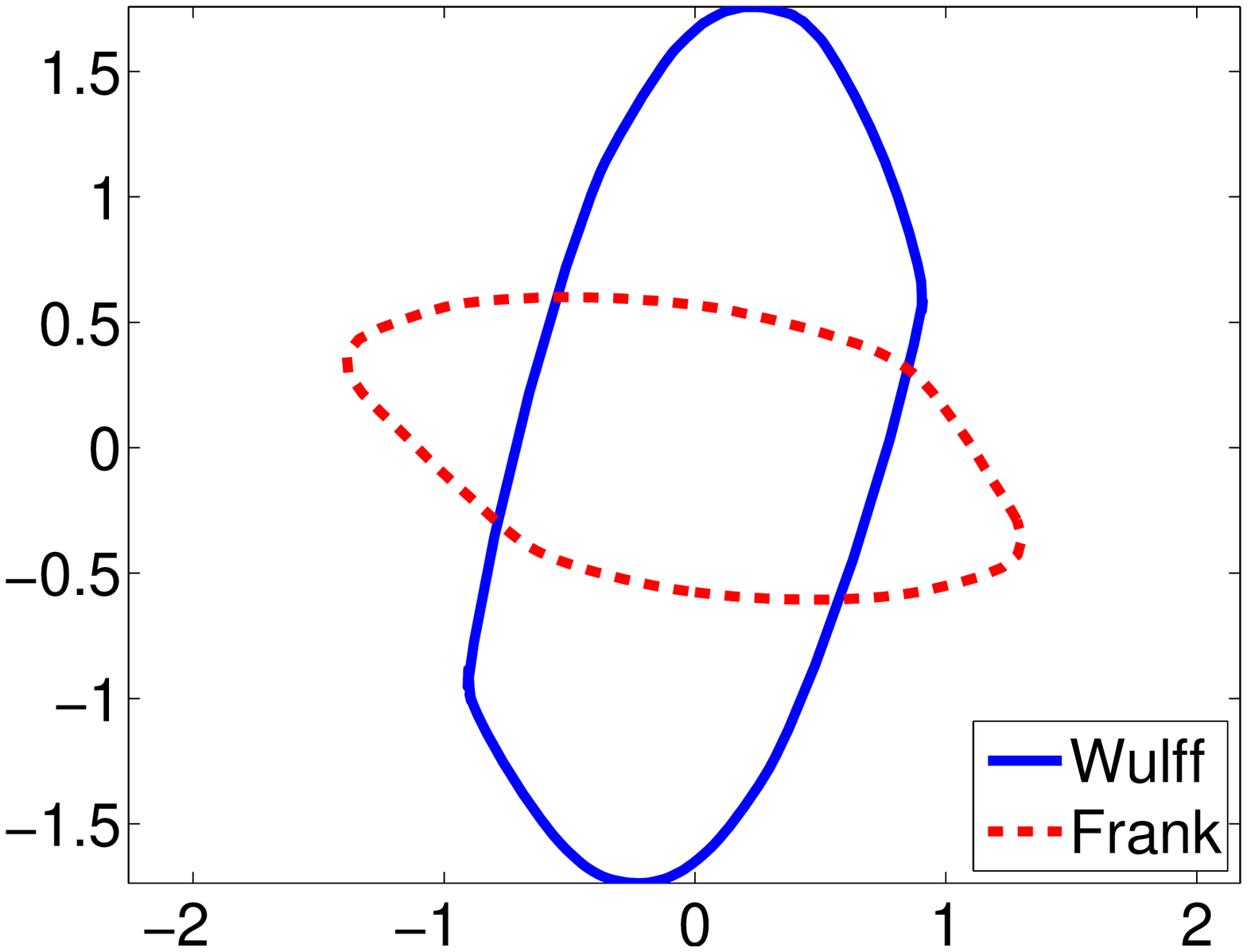}}
\end{center}
\caption{\small
(a) the curve $\Gamma$; (b) the optimal anisotropy function $\sigma\equiv\sigma^N, N=50$; (c) the Wulff shape $W_\sigma$ and Frank diagram ${\mathcal F}_\sigma$ 
}
\label{fig:batman}
\end{figure}

In Figure~\ref{fig:batman} (a) we present a simple test example of a Jordan curve $\Gamma=\{\vecx(u) \ | \ u\in[0,1]\}$ where $\vecx(u)=(x_1(u), x_2(u))^T$, and
\[
x_1(u) = \cos(2\pi u), \ 
x_2(u) = 0.7\sin(2\pi u) + \sin(\cos(2\pi u)) + \left(\sin(6\pi u)\sin(2\pi u)\right)^2.
\]
The curve was discretized by $K=1000$ grid points and the Fourier length spectrum coefficients were computed according to \eqref{cp-discr}. We chose $N=50$ Fourier modes in this example. The anisoperimetric ratio for the optimal anisotropy function $\sigma$ (depicted in Figure~\ref{fig:batman} (b) equals 2.306 whereas the isoperimetric ratio of $\Gamma$ equals 3.041. The Wulff and Frank diagrams are shown in In Figure~\ref{fig:batman} (c).
In Table~\ref{tab1} we present results of computation for various numbers $N$ of Fourier modes for a curve shown in Fig~\ref{fig:batman} (a). The area $|W_{\sigma^N}|$ of the optimal Wulff shape converges to the value 4.1612 as $N\approx 300$ when we impose the constraint $L_{\sigma^N}(\Gamma)=L(\Gamma)=9.167$. It should be also noted that satisfactory numerical results were obtained for rather low dimensions $N\approx 50$.

It is known that the worst case time complexity of SeDuMi implementation (including main and inner iterations) is $O(n_v^2 n_c^{2.5}+ n_c^{3.5})$ where $n_c$ and $n_v$ are the numbers of variables and constraints, respectively (c.f. \cite{henrion}). Since the number of constraints $n_c=O(N)$ and number of variables $n_v=O(N^2)$ (see Table~\ref{tab1}) the worst case time complexity should have the order $O(N^{6.5})$.  We calculated the experimental order of time complexity (eotc) by comparing elapsed times $T_k$ for different $N_k$ as follows: $eotc_k = \ln(T_{k+1}/T_k)/\ln(N_{k+1}/N_k)$. It turns out the $eotc\lesssim 3.8$, i.~e. $T \lesssim O(N^{3.8})$, so it is below the worst case complexity. All
computations were performed on Quad-Core AMD Opteron Processor with 2.4~GHz frequency, 32~GB of memory. 

\begin{figure}
\begin{center}
\subfloat[]{\includegraphics[width=0.28\textwidth]{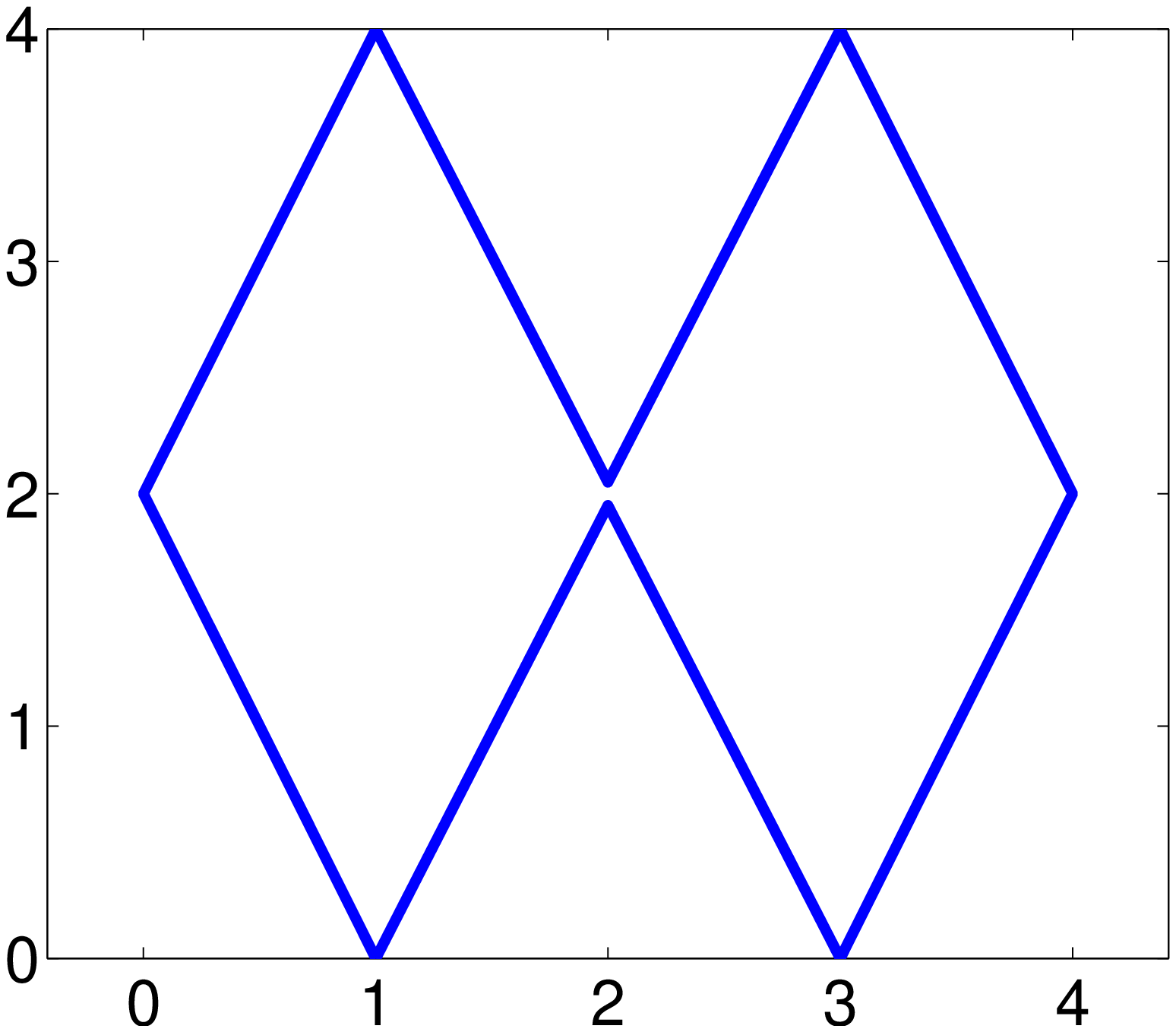}}
\subfloat[]{\includegraphics[width=0.28\textwidth]{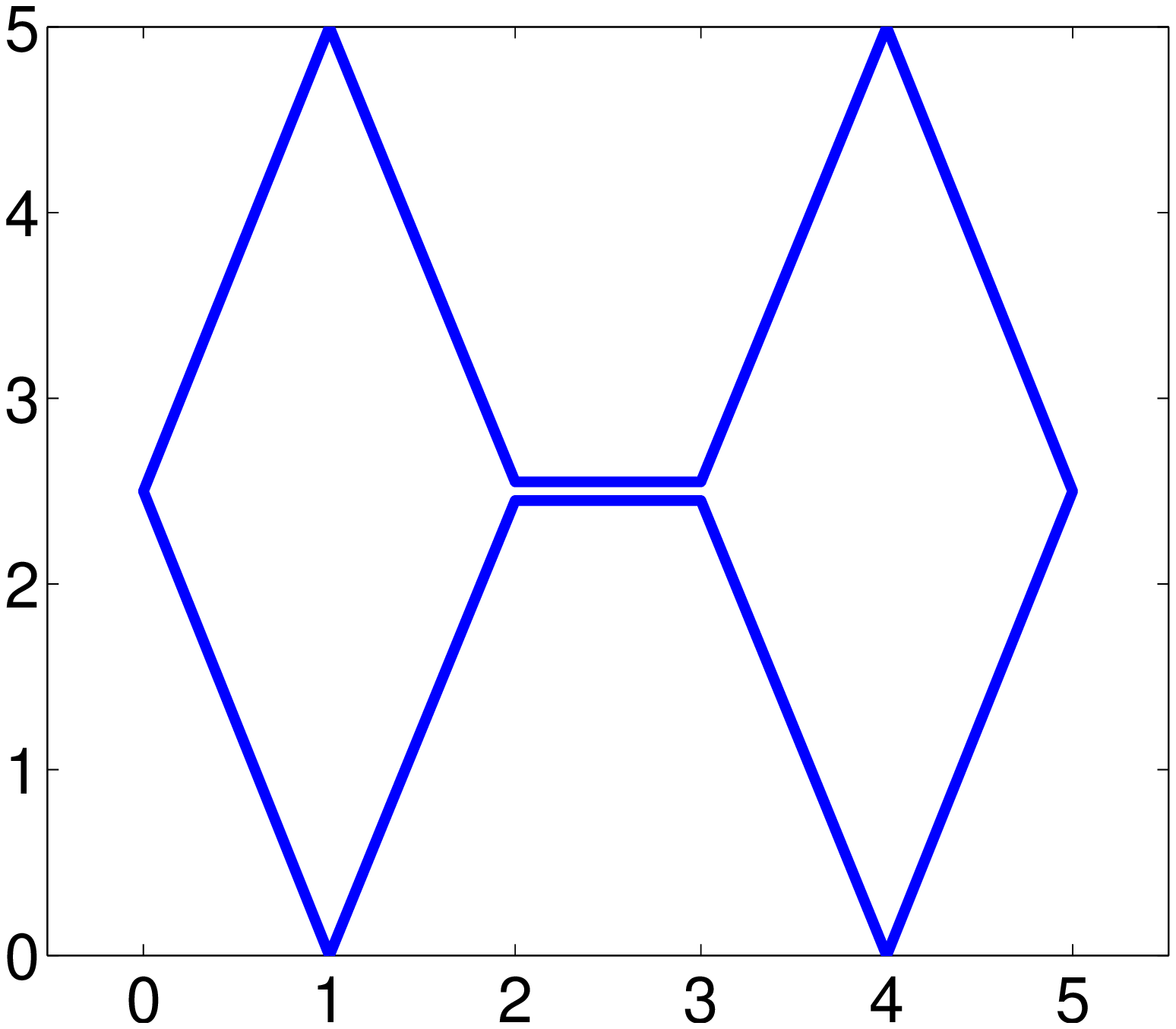}}
\subfloat[]{\includegraphics[width=0.28\textwidth]{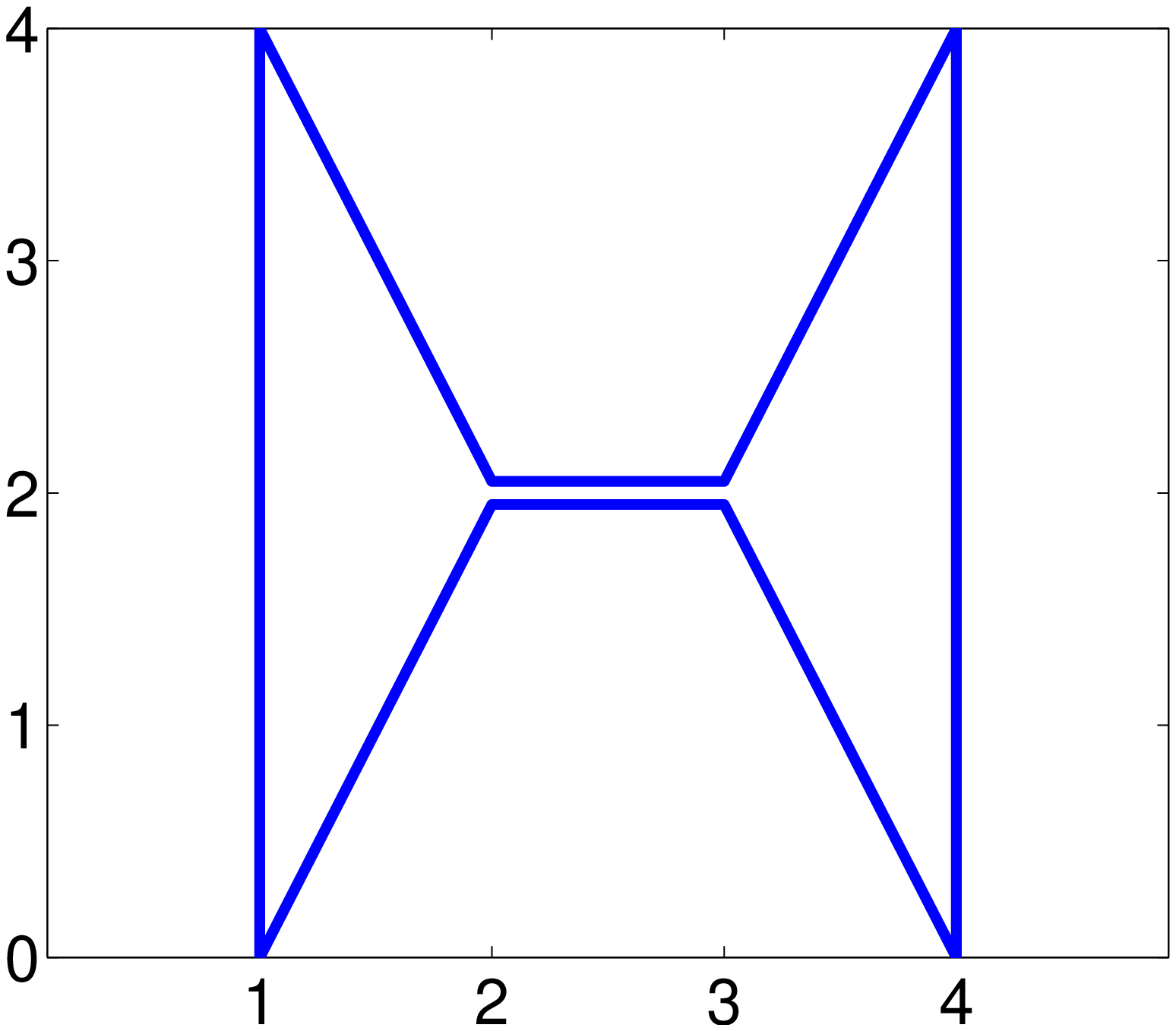}}

\subfloat[]{\includegraphics[width=0.28\textwidth]{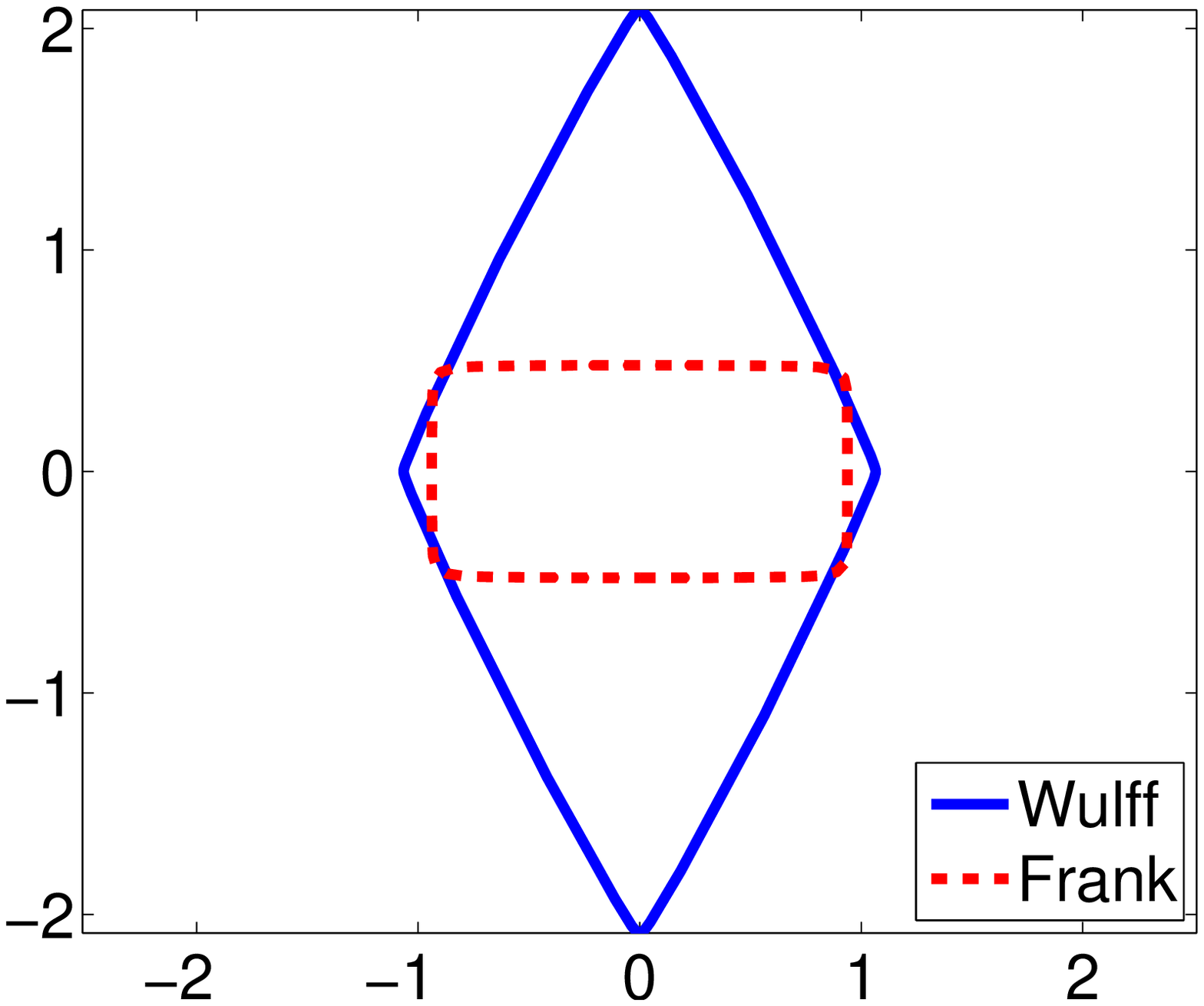}}
\subfloat[]{\includegraphics[width=0.28\textwidth]{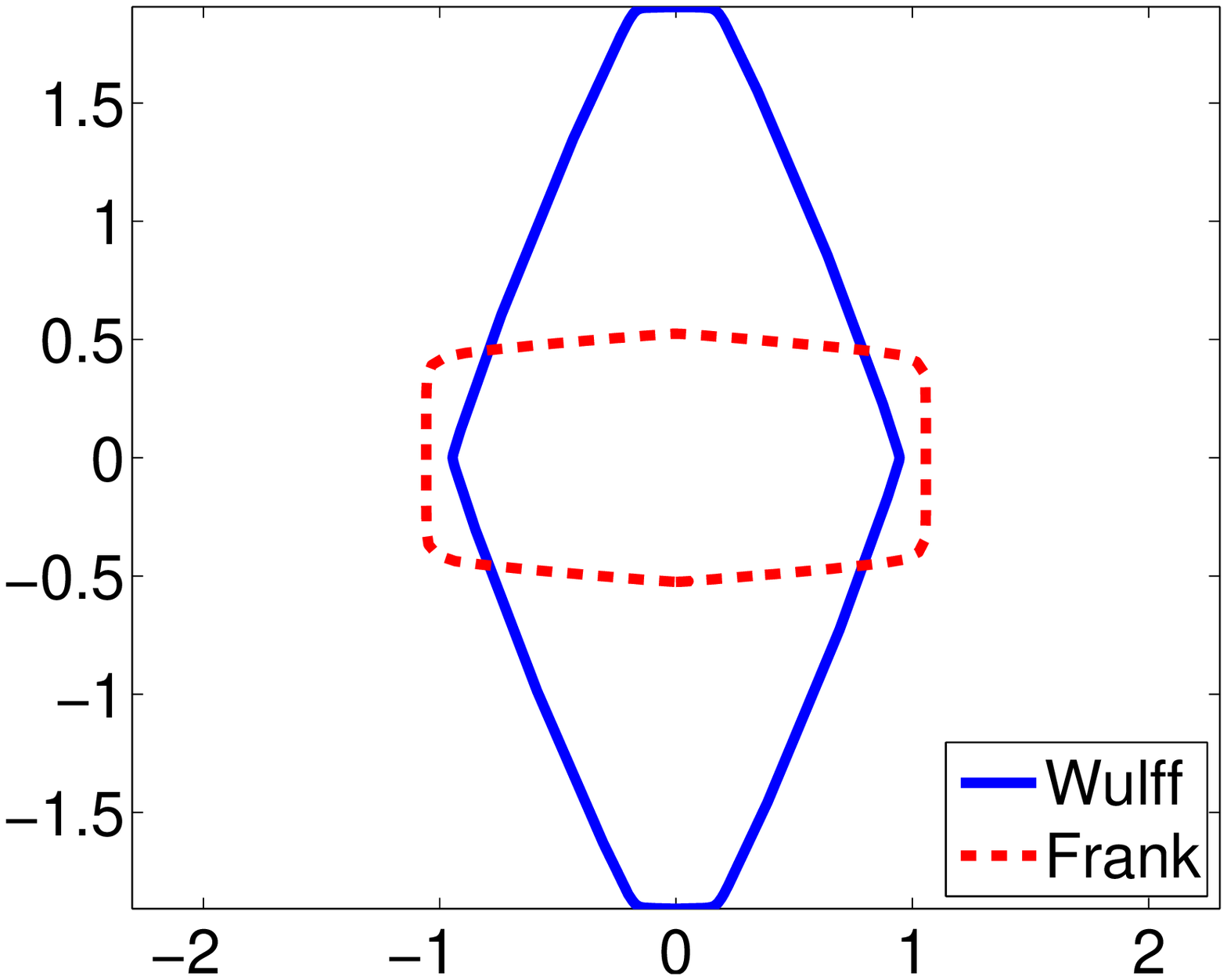}}
\subfloat[]{\includegraphics[width=0.28\textwidth]{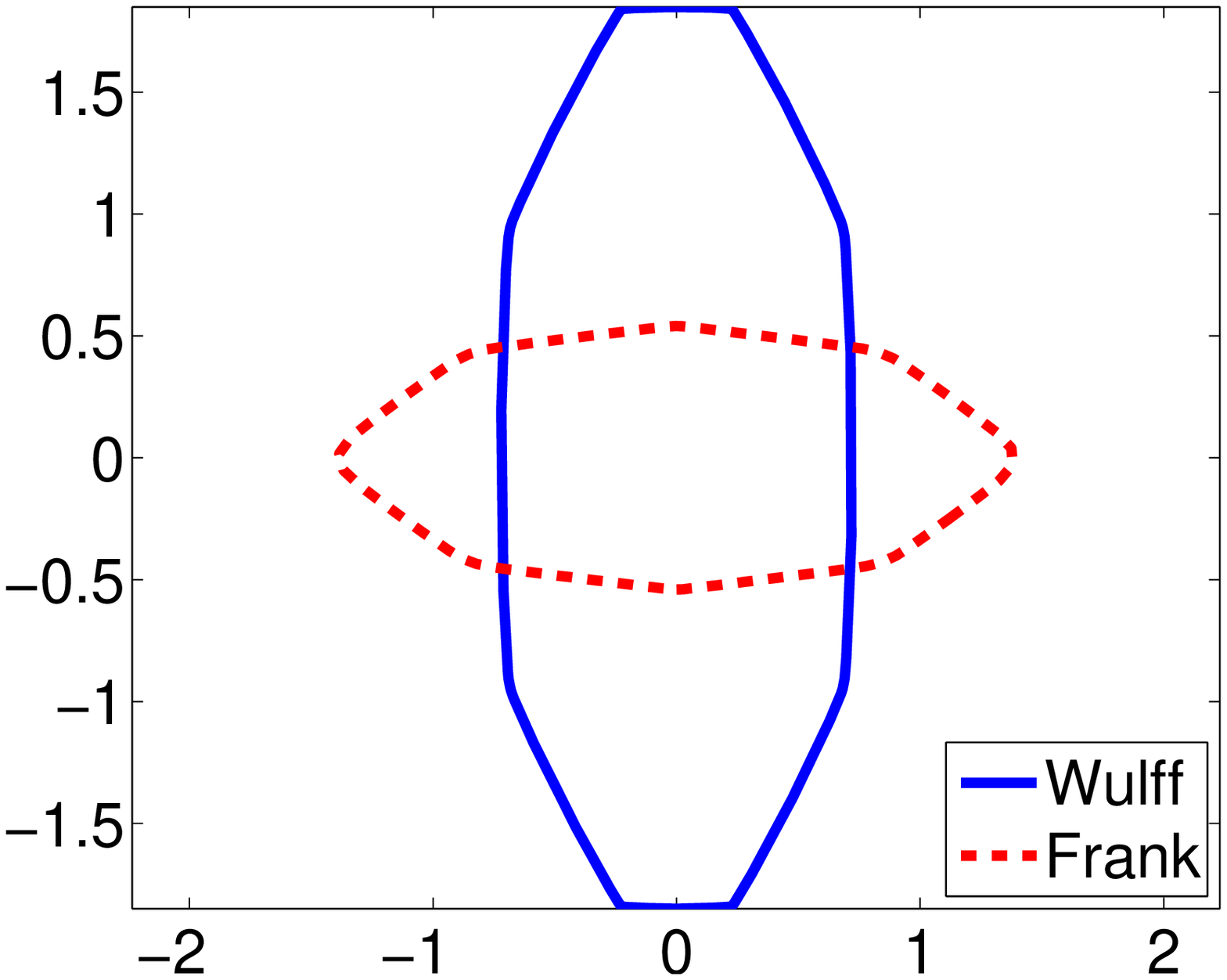}}
\end{center}
\caption{\small
Nonconvex polygonal curves (a-c) and their corresponding optimal Wulff shapes and Frank diagrams (d-f)
}
\label{fig:polygon}
\end{figure}

In Figure~\ref{fig:polygon} we present results of resolution of the optimal anisotropy function for various polygonal curves (a-c). The corresponding optimal Wulff shapes and Frank diagrams (e-f) show their anisotropy structure. For instance, there are four outer normal directions of facets in Figure~\ref{fig:polygon} (a). The corresponding Wulff shape in Figure~\ref{fig:polygon} (d), solid blue line, has a shape of the four fold anisotropy with the same set of outer normal directions. Similarly, other polygons shown in Figure~\ref{fig:polygon} have hexagonal (b-e) and octagonal (c-f) anisotropy and the sets of their outer normal vectors to facets coincide. We again chose a sufficiently large number $N=50$ of Fourier modes in these examples so that numerically computed Wulff shapes are just slightly rounded polygons. 

\begin{figure}
\begin{center}
\subfloat[]{\includegraphics[width=0.24\textwidth]{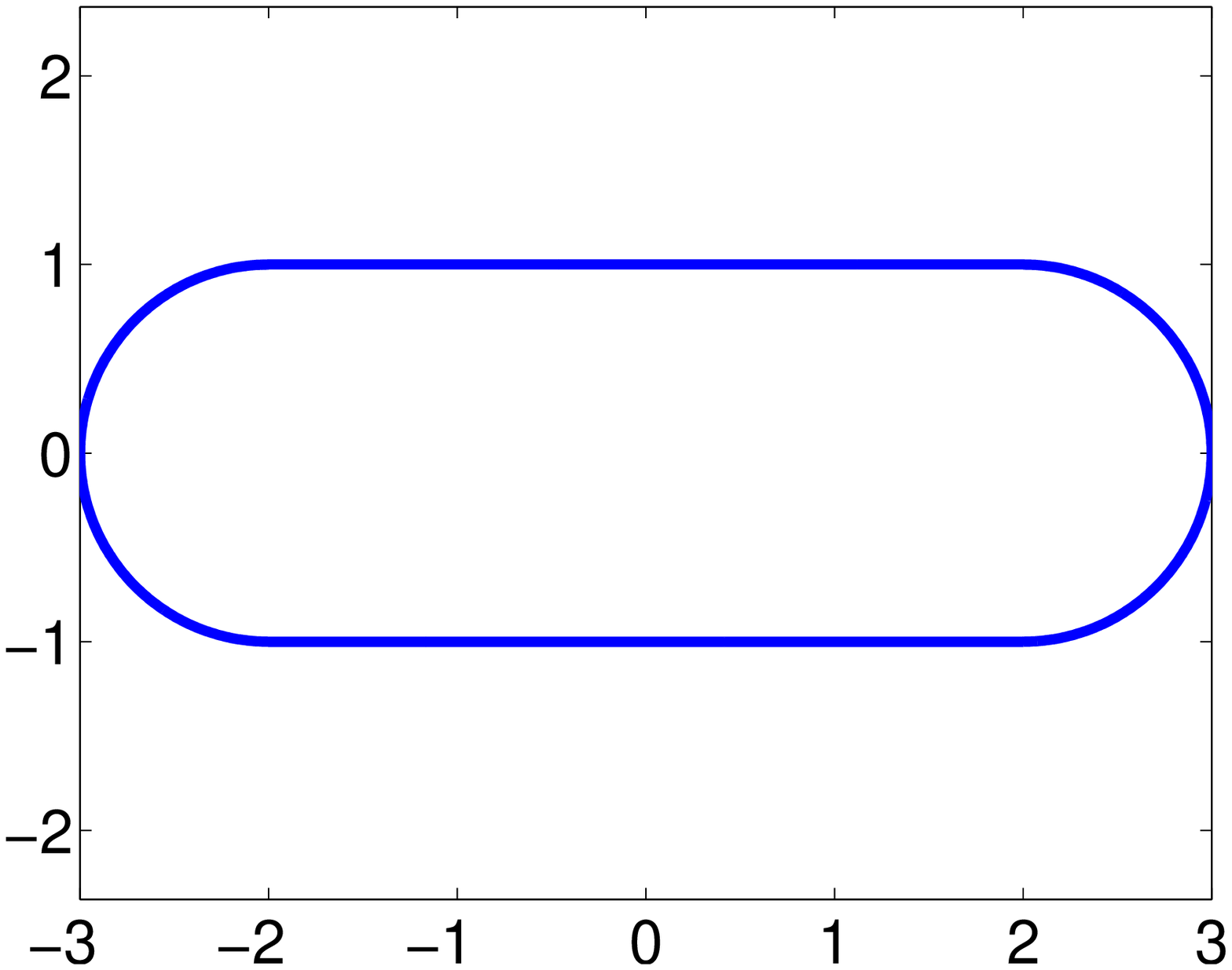}}
\subfloat[]{\includegraphics[width=0.24\textwidth]{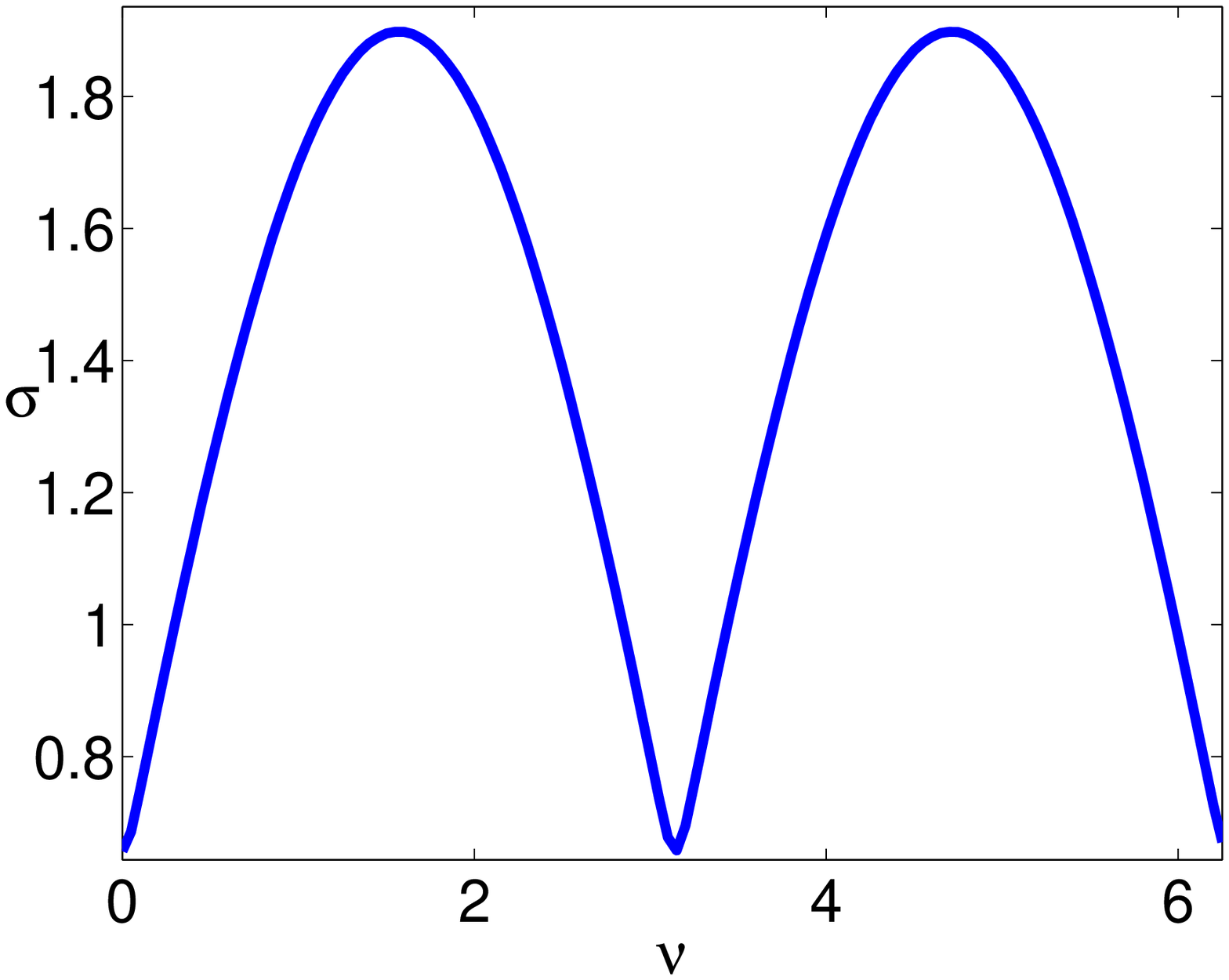}}
\subfloat[]{\includegraphics[width=0.24\textwidth]{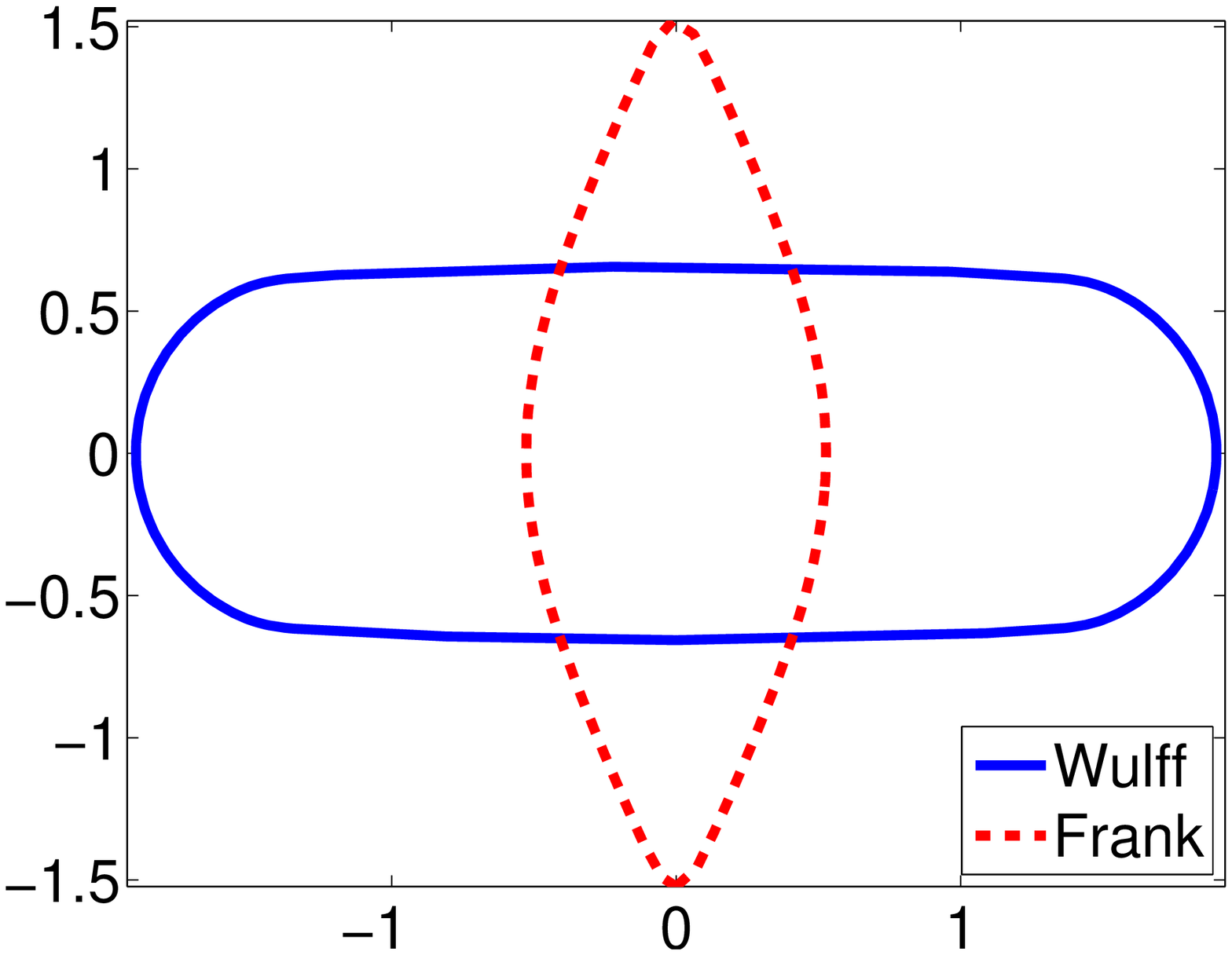}}
\subfloat[]{\includegraphics[width=0.24\textwidth]{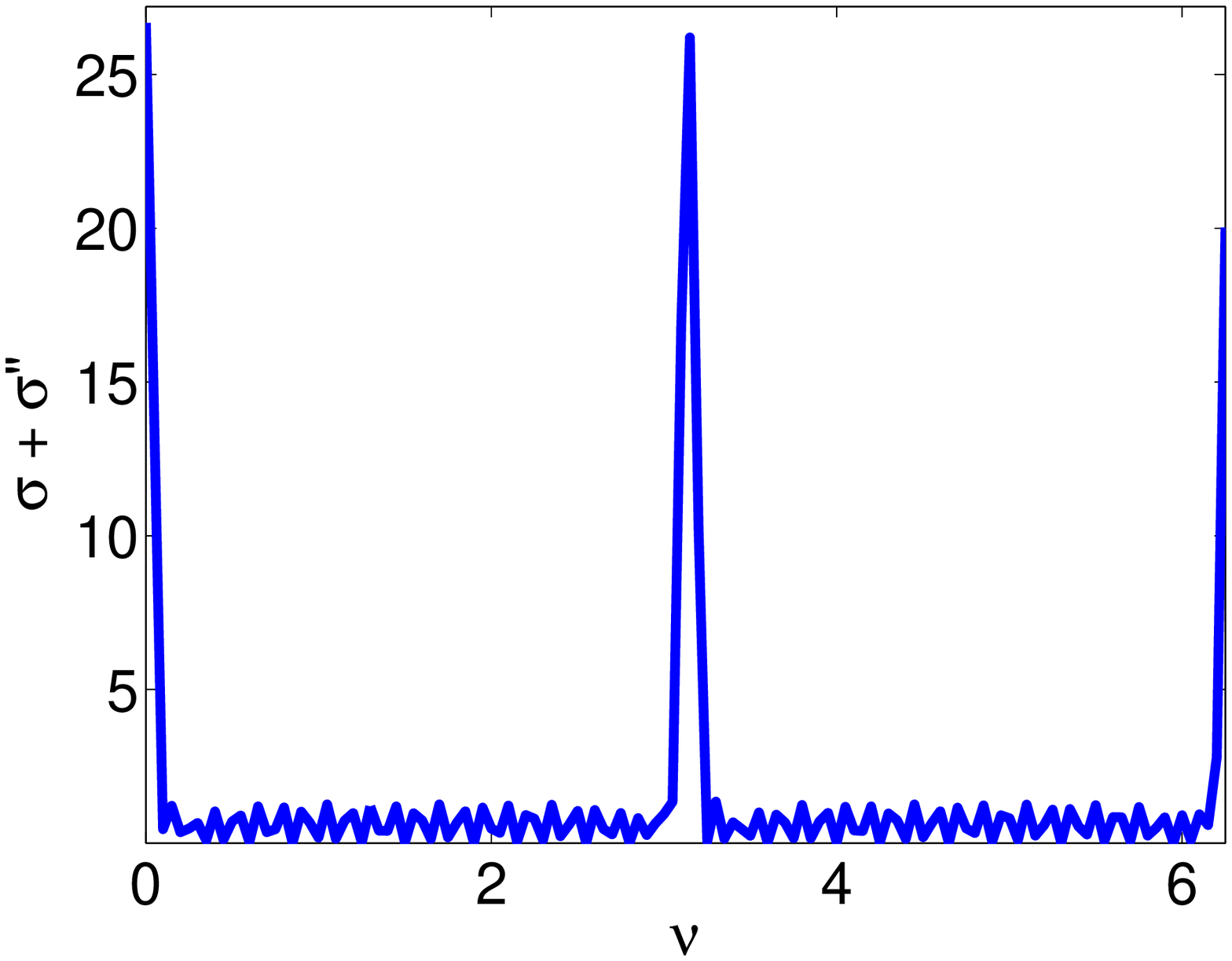}}
\end{center}
\caption{\small
(a) the "capsule" curve $\Gamma$; (b) the optimal anisotropy function $\sigma\equiv\sigma^N, N=50$; (c) the Wulff shape $W_\sigma$ and Frank diagram ${\mathcal F}_\sigma$;  
(d) the graph of the function $\sigma+\sigma''\equiv 1/\kappa$
}
\label{fig:capsule}
\end{figure}

\begin{table}[ht]
\caption{Dependence of the time complexity of computation with respect to the number of Fourier modes $N$, its experimental order of time complexity (eotc) and the area of the  optimal Wulff shape $|W_{\sigma^N}|$}
\label{tab1}
\scriptsize
\begin{center}
\begin{tabular}
{lllllllll}
\hline\noalign{\smallskip}
$N$& 
$n_c$& 
$n_v$& 
time (s)& 
eotc& 
$|W_{\sigma^N}|$ \\
\noalign{\smallskip}\hline\noalign{\smallskip}
25& 
102& 
2576& 
1& 
--& 
4.14087 \\
50& 
202& 
10151& 
6& 
2.58& 
4.14519 \\
100& 
402& 
40301& 
47& 
2.97& 
4.14597 \\
200& 
802& 
160601& 
434& 
3.21& 
4.14611 \\
300& 
1202& 
360901& 
1692& 
3.36& 
4.14612 \\
350& 
1402&
491051&
3020& 
3.75& 
4.14624 \\
\noalign{\smallskip}\hline
\end{tabular}

\end{center}

\end{table}

For any strictly convex $C^2$ curve $\Gamma$ the optimal anisotropy function $\sigma$ corresponds to the Wulff shape with $\partial W_\sigma = \Gamma$ and $\Pi_\sigma(\Gamma)=1$. On the other hand, if $\Gamma$ is just a piecewise $C^2$ smooth curve then there need not exist a minimizing anisotropy function belonging to $\K\subset W^{2,2}_{per}(0,2\pi)$. The purpose of the next example shown in Figure~\ref{fig:capsule} is to illustrate behavior of Sobolev norms in the space $W^{k,2}_{per}(0,2\pi), k=0,1,2,$ of the optimal solution $\sigma^N$ for various dimensions $N\in\N$. We consider the "capsule" like curve $\Gamma$ from Example~\ref{fourier-spec-ex-capsule} with $l=4, r=1$. This a $C^1$ smooth and only piecewise $C^2$ smooth convex curve. According to Remark~\ref{rem-capsule} we have $\inf_{\sigma\in\K} \Pi_\sigma(\Gamma)=1$ but there is no minimizer $\sigma$ belonging to $\K\subset W^{2,2}_{per}(0,2\pi)$. It can be deduced from Table~\ref{tab2} that the Sobolev norms $\Vert\sigma^N\Vert_{k,2}, k=0,1,$ stay bounded for $N\gg 1$. On the other hand, the norm $\Vert\sigma^N\Vert_{2,2}$ becomes unbounded and $\Vert \sigma^N \Vert_{2,2} \approx O(N^{0.5})$ as $N\to\infty$, i.~e. the experimental order of convergence is approximately $0.5$. The pointwise behavior of the function $\sigma+\sigma''\equiv1/\kappa$ is shown in  Figure~\ref{fig:capsule}, (d). The anisoperimetric ratio $\Pi_{\sigma^N}(\Gamma)$ tends to unity with the speed of $O(N^{-1})$, i.~e. the experimental order of convergence is $-1$.

\begin{table}[ht]
\caption{Dependence of various Sobolev norms of the optimal anisotropy function with respect to the dimension $N$. Experimental rate of divergence of the $\Vert.\Vert_{2,2}$ norm and experimental rate of convergence of the anisoperimetric ratio to unity
}
\label{tab2}
\scriptsize
\begin{center}
\begin{tabular}
{lllllll}
\hline\noalign{\smallskip}
$N$& 
$\Vert \sigma^N \Vert_{0,2}$&
$\Vert \sigma^N \Vert_{1,2}$&
$\Vert \sigma^N \Vert_{2,2}$& 

eoc($W^{2,2}_{per}$) & 
$\Pi_{\sigma_N} -1 $& 
eoc($\Pi_\sigma -1$) 
\\
\noalign{\smallskip}\hline\noalign{\smallskip}
10 & 1.3937 & 1.5032 & 2.0164 & -- & 0.057060 & -- \\
25 & 1.4467 & 1.5701 & 2.6417 & 0.29 & 0.023216 & -0.98 \\
50 & 1.4634 & 1.5910 & 3.2631 & 0.3 & 0.012790 & -0.86 \\
100 & 1.4730 & 1.6029 & 4.2039 & 0.37 & 0.006784 & -0.91 \\
150 & 1.4763 & 1.6071 & 4.9362 & 0.4 & 0.004665 & -0.92 \\
200 & 1.4780 & 1.6093 & 5.5618 & 0.41 & 0.003567 & -0.93 \\
250 & 1.4791 & 1.6106 & 6.1141 & 0.42 & 0.002894 & -0.94 \\
\noalign{\smallskip}\hline
\end{tabular}

\end{center}

\end{table}

\begin{figure}
\begin{center}

\subfloat[]{\includegraphics[width=0.31\textwidth]{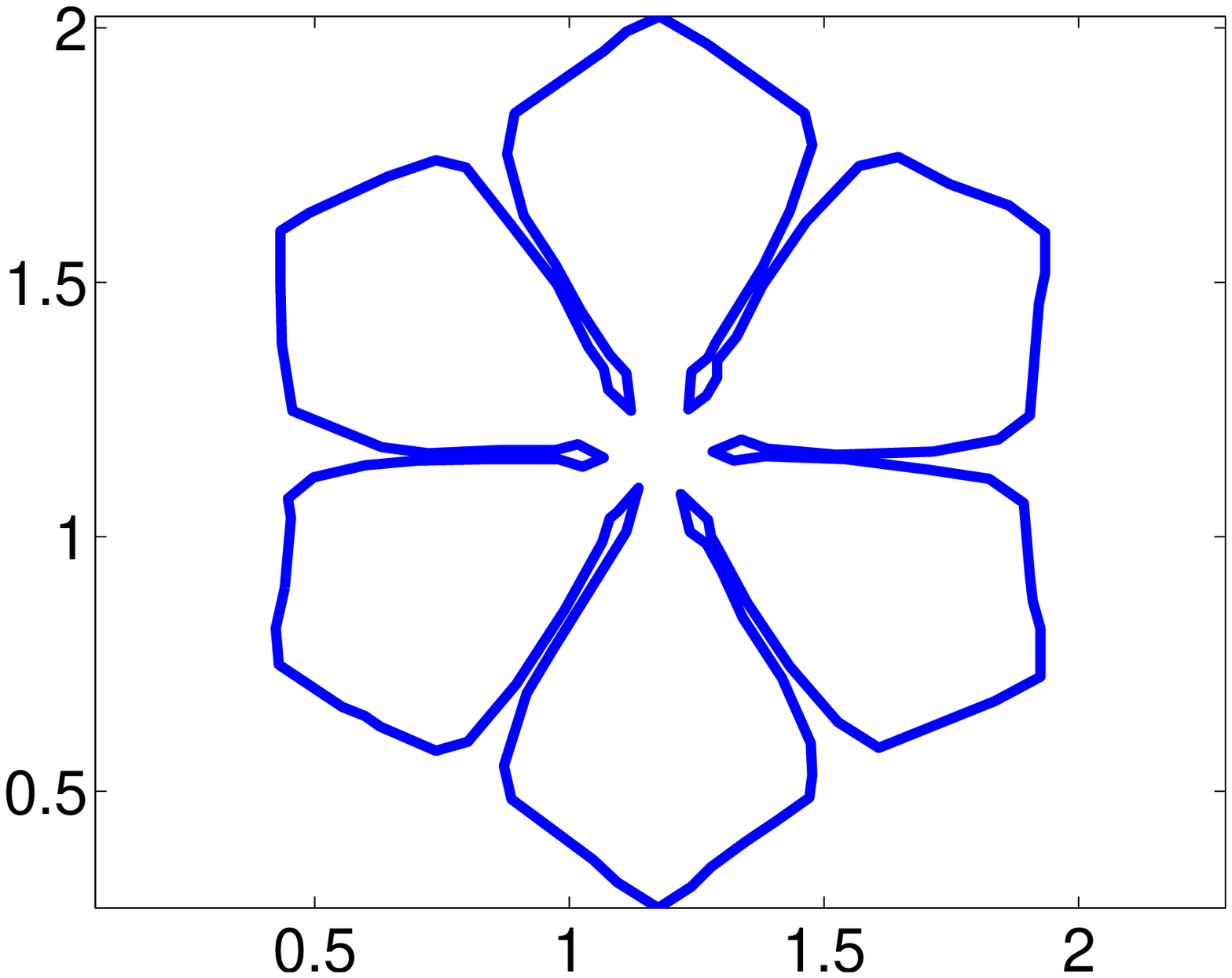}}
\subfloat[]{\includegraphics[width=0.31\textwidth]{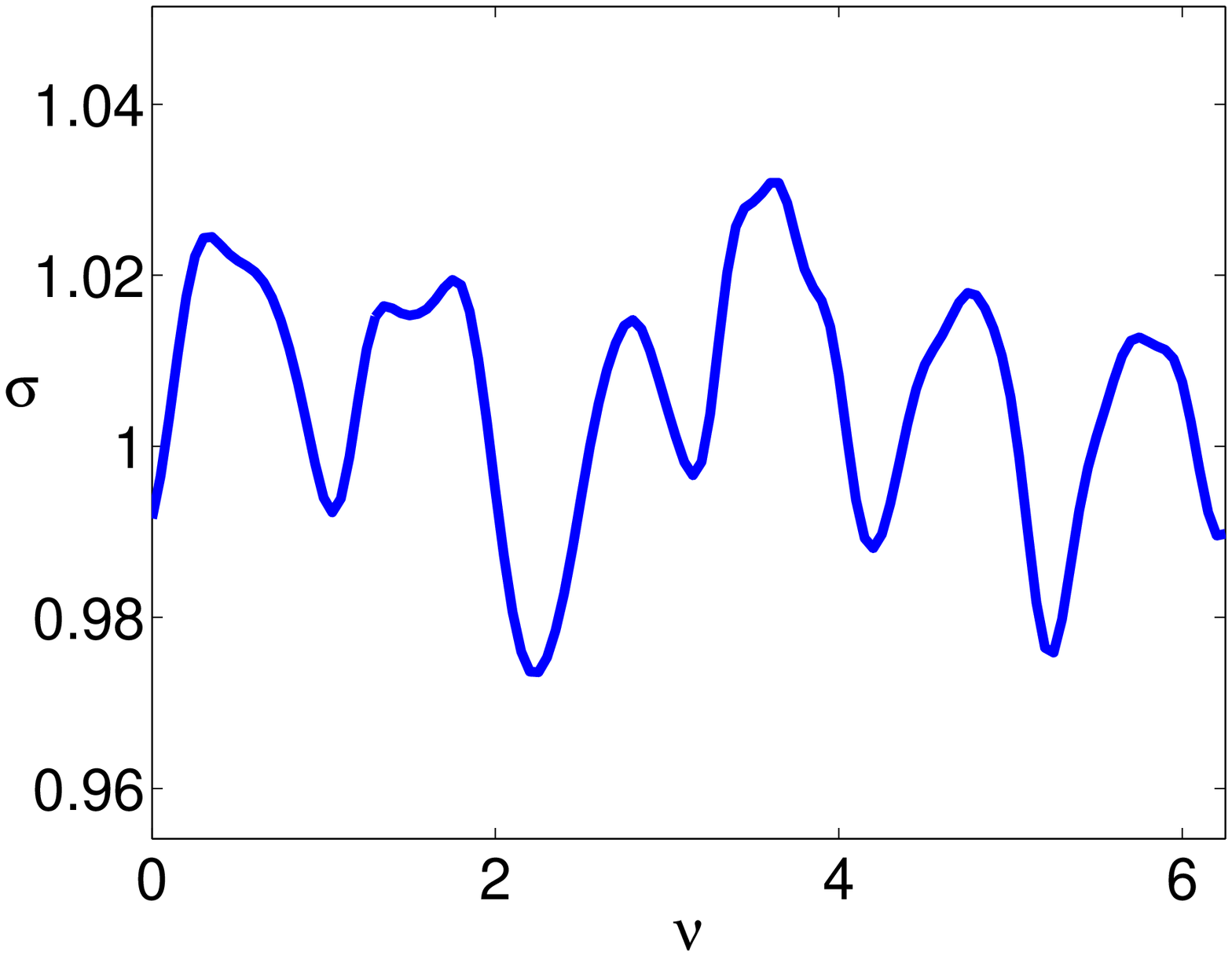}}
\subfloat[]{\includegraphics[width=0.31\textwidth]{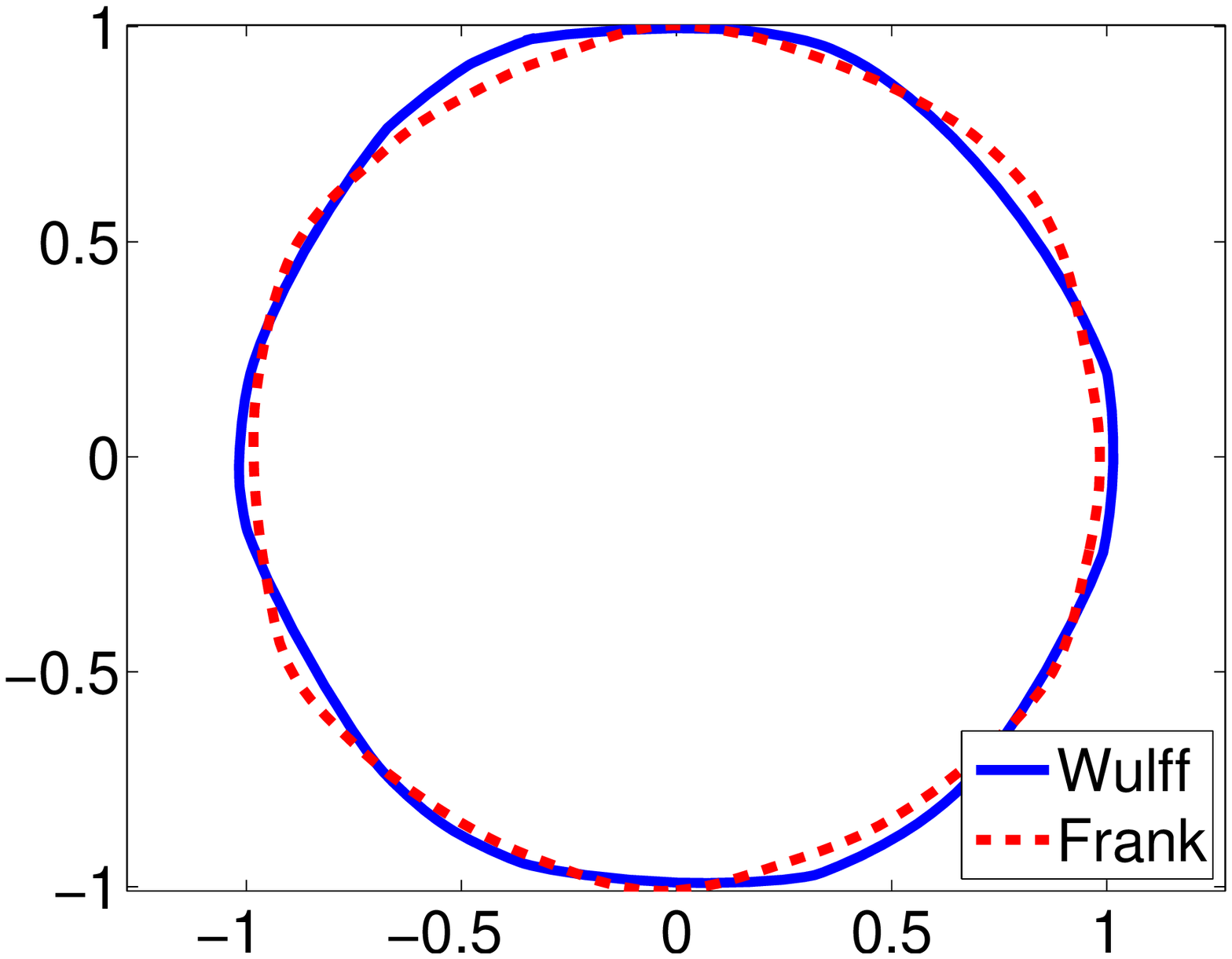}}
\\
\subfloat[]{\includegraphics[width=0.31\textwidth]{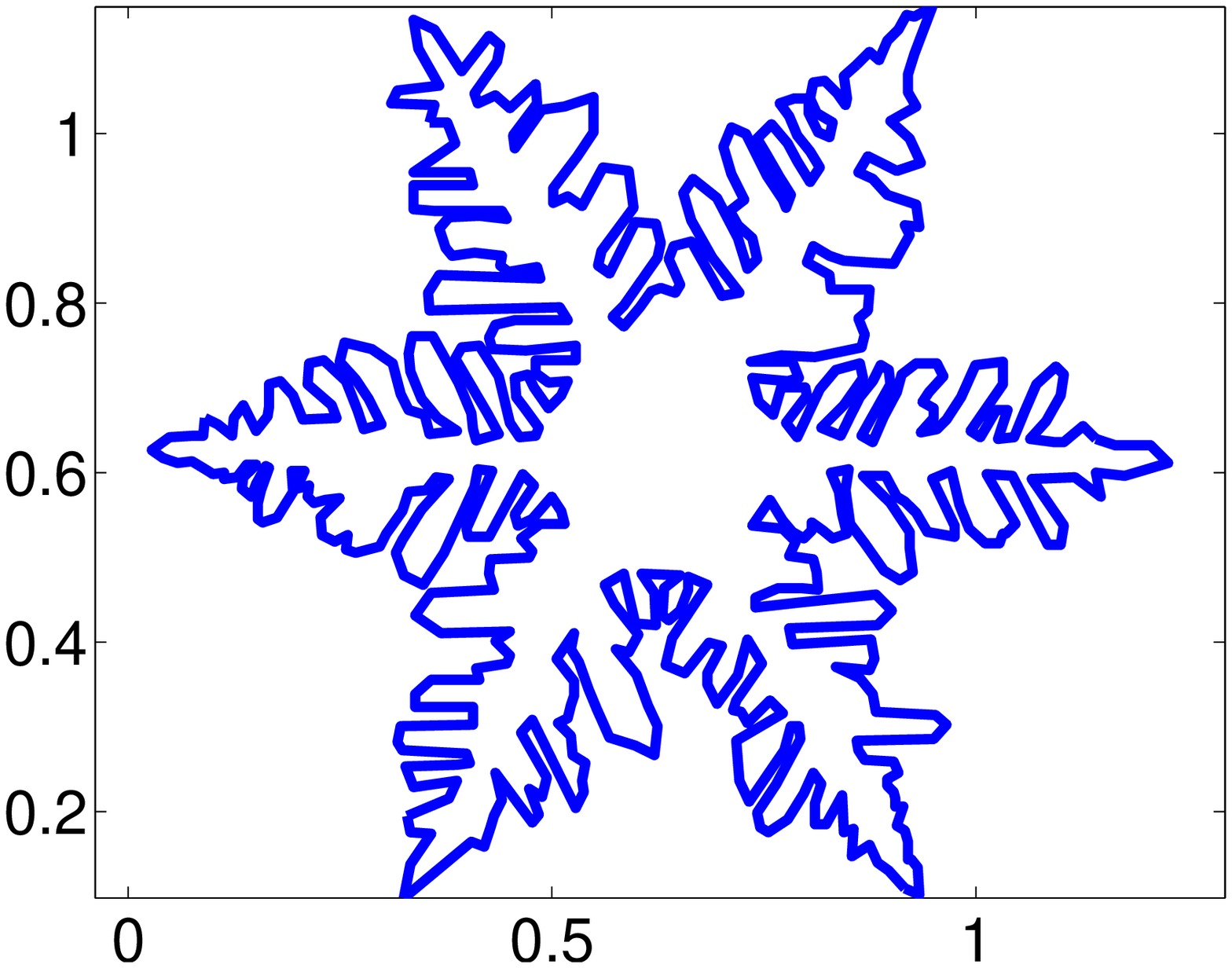}}
\subfloat[]{\includegraphics[width=0.31\textwidth]{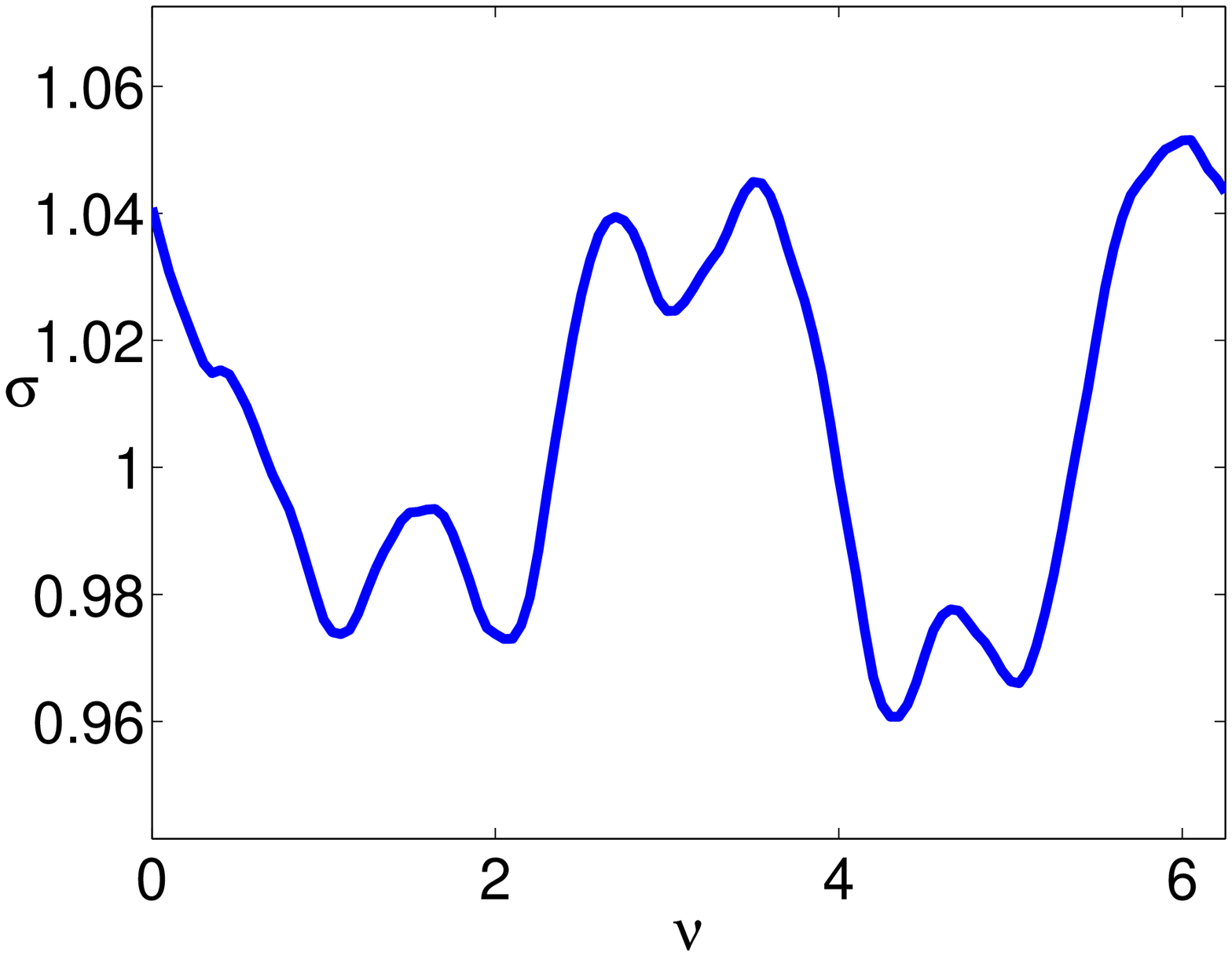}}
\subfloat[]{\includegraphics[width=0.31\textwidth]{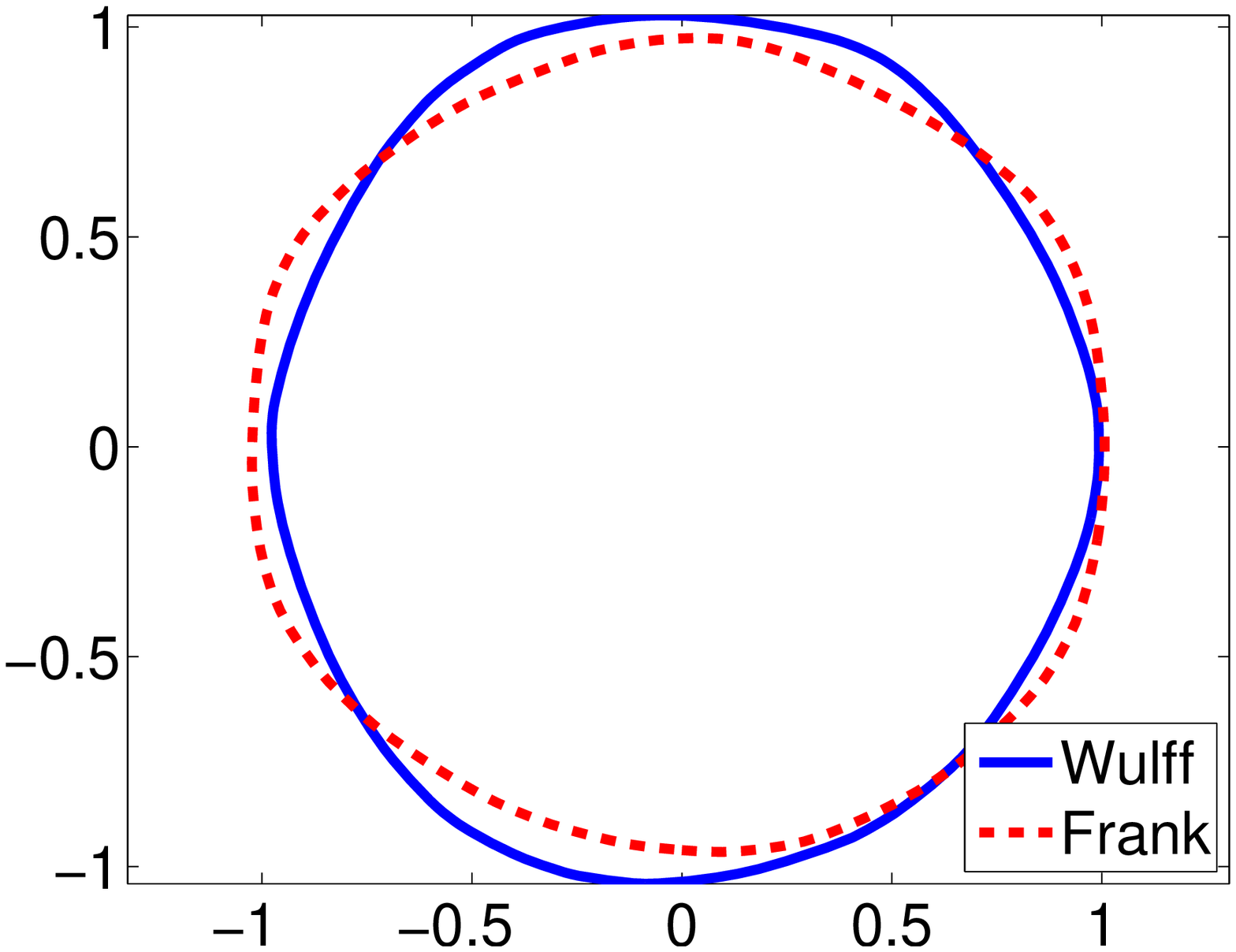}}
\end{center}
\caption{\small
Jordan curves (a,d) representing to boundaries of snowflakes. The optimal anisotropy function $\sigma^N$ is shown in (b,e). The Wulff shape and Frank diagram are depicted in (c,f)
}
\label{fig:snowflake6}
\end{figure}

In the last two numerical examples shown in Figure~\ref{fig:snowflake6} we present computation of the optimal anisotropy functions for boundaries of real snowflakes. We used $N=50$ Fourier modes. In both cases we resolved the optimal anisotropy function $\sigma$ corresponding to hexagonal symmetry, as it can be expected for snowflake crystal growth. If we introduce the anisotropy strength as follows:
\[
\varepsilon := (\sigma_{max} - \sigma_{min})/(2\sigma_{avg}),
\]
where $\sigma_{max}, \sigma_{min}$ are the maximal and minimal values of $\sigma(\nu)$ for $\nu\in[0,2\pi]$ and $\sigma_{avg} = \frac{1}{2\pi}\int_0^{2\pi} \sigma(\nu)\dnu = \sigma_0$ then $\varepsilon=0.028850$ for the snowflake (a) and $\varepsilon=0.045666$ for the snowflake (d)\footnote{Snowflake images sources:\\
(a) http://milanturek.files.wordpress.com/2010/07/sneh-vlocka8.jpg \\
(d) http://www.isifa.com/data/dispatches/ed/167/\_main.jpg
}. Notice that the maximal value of the strength parameter $\varepsilon$ for the anisotropy function of the form: $\sigma(\nu)=1+\varepsilon\cos(m\nu)$ is $\varepsilon=1/(m^2-1) = 0.0285714$ for hexagonal symmetry $m=6$.

\section{Conclusions}
We proposed a new method for resolving the optimal anisotropy function that minimizes the anisoperimetric ratio for a given Jordan curve in the plane. Construction of the optimal anisotropy function can be regarded as a solution to the inverse Wulff problem. Our approach of solving the inverse Wulff problem was based on reformulation of the optimization problem in terms of complex Fourier coefficients of the anisotropy functions. We furthermore proposed and analyzed the Fourier length spectrum of a curve. Using results from the theory of semidefinite matrices we were able to prove useful asymptotic estimates on elements of the Fourier length spectrum. It turned out that the finite Fourier modes approximation leads to an indefinite quadratic optimization problem with linear matrix inequalities. We solved this problem by means of the so-called enhanced semidefinite relaxation method. It consisted in solving the relaxed convex semidefinite problem obtained as the second Lagrangian dual of the original problem augmented by a quadratic-linear constraint. Various numerical examples and tests of experimental order of convergence were presented. In particular, we presented examples of computation of optimal anisotropy function for a set of snowflake boundaries. 

\section*{Acknowledgments}
The authors thank Prof. M. Halick\'a and M. Hamala for their constructive comments and suggestions concerning semidefinite relaxation methods.

\section*{Appendix}

In this appendix section we provide a detailed derivation of the first and second Lagrangian duals of \eqref{P1} and \eqref{D2}, respectively. 

Clearly, $\hbox{tr}(V^TAxx^T)=x^TV^TAx=x^TA^TVx=\frac12x^T(V^TA+A^TV)x$. Since $\hbox{tr}(V^Tbx^T)=x^TV^Tb=b^TVx$ the Lagrangian ${\mathcal L}^1$ can be rewritten in a compact form as follows:
${\mathcal L}^1(x;\lambda, u, V, Z )=x^TQx+2s^Tx+\tau$, where 
$Q=P_0+\sum_{l=1}^d\lambda_l P_l+\frac12(V^TA+A^TV), \ s=q_0+ \sum_{l=1}^d\lambda_l q_l  +\frac12(A^Tu-V^Tb-z)$, 
$\tau=r_0+ \lambda^T r - u^T b - z_0,$ with $z=(z_1,\ldots, z_n), \ z_j=\hbox{tr}(Z^T\tilde{H}_j)$ and $r=(r_1, \ldots, r_l)^T$. 

Then the Lagrange dual function $\mathcal{G}^1$ can be defined as: \\ $\mathcal{G}^1(\lambda, u, V, Z)=\inf_x{\mathcal L}^1(x; \lambda, u, V, Z)$. Since $\partial_x {\mathcal L}^1=2Qx+2s,$ the dual function attains a finite value $\mathcal{G}^1 >-\infty$ if and only if $Q$ is positive semidefinite and the vector $s$ belongs to the range $\mathcal{R}[Q]$ of $Q$. If these two conditions are satisfied then $\mathcal{G}^1=-s^TQ^{\dagger} s +\tau$ where $Q^{\dagger}$ is the Moore-Penrose pseudoinverse to $Q$. The dual problem has the form:
\[
\begin{array}{rl}
\max & \mathcal{G}^1(\lambda, u, V, Z)\\
\hbox{s. t. } & q_0+\lambda q_1+\frac12(A^Tu-V^Tb-z)\in \mathcal{R}[P_0+ \sum_{l=1}^d\lambda_l P_l+\frac12(V^TA+A^TV)],\\
  & P_0+\sum_{l=1}^d\lambda_l P_l+\frac12(V^TA+A^TV)\succeq 0,\\
  & Z\succeq 0, \lambda\geq 0,\ \ z_j=\hbox{tr}(Z^T\tilde{H}_j),\  j=0,1,\ldots, n.
\end{array}
\]
Using the properties of the generalized Schur complement (c.f. \cite{zhang}), it can be rewritten in the following form:
\begin{equation*}
\begin{array}{rl}
\max & \gamma\\
\hbox{s. t.} 
 & Z\succeq 0, \lambda\geq 0,\ z_j=\hbox{tr}(Z^T\tilde{H}_j),\  j=0,1,\ldots, n, \\
  & \hskip-1.2truecm \left(\begin{array}{cc}
P_0+ \sum_{l=1}^d \lambda_l P_l+\frac12(V^TA+A^TV), & q_0+\lambda^T q +\frac12(A^Tu-V^Tb-z)\\
(q_0+ \lambda^T q + \frac12(A^Tu-V^Tb-z))^T, & r_0+\lambda^T r -u^Tb-z_0-\gamma\\
\end{array}\right)\succeq 0.\\
\end{array}
\end{equation*}
If we introduce the notation for $M_0, M_l, N_0, N_j$ and $M_*$ then we obtain \eqref{D2}.

Next we derive the second dual problem. To this end, consider the  Lagrangian ${\mathcal L}^2(\gamma,\lambda, Z,V,u,z; W,\beta,\tilde{X},\alpha)$  for problem \eqref{D2}, i.~e.
\begin{eqnarray*}
{\mathcal L}^2 &=&\gamma+\hbox{tr}(ZW)+\lambda\beta
+\sum_{j=0}^n\alpha_j(z_j-\hbox{tr}(Z\tilde{H}_j))
\\
&&+\hbox{tr}(\tilde{X}(M_0+ \sum_{l=1}^d \lambda_l M_l + M_*(u,V) -\sum_{j=0}^n z_jN_j-\gamma N_0))
\\
&=&\gamma(1-\hbox{tr}(\tilde{X}N_0))+\sum_{l=1}^d \lambda_l[\beta_l+\hbox{tr}(\tilde{X}M_l)] + \hbox{tr}(Z(W-\sum_{j=0}^n\alpha_j\tilde{H}_j))
\\
&&
+\hbox{tr}(\tilde{X} M_*(u,V)) + \hbox{tr}(\tilde{X}M_0)+\sum_{j=0}^n z_j(\alpha_j-\hbox{tr}(\tilde{X}N_j)).
\end{eqnarray*}
The function ${\mathcal L}^2(\gamma,\lambda, Z,V,u,z; W,\beta,\tilde{X},\alpha)$ is linear in the  $(\gamma,\lambda, Z,V,u,z)$ variable The dual function is defined as follows:
\[
\mathcal{G}^2( W,\beta,\tilde{X},\alpha)=\sup_{\gamma,\lambda, Z,V,u,z} {\mathcal L}^2(\gamma,\lambda, Z,V,u,z; W,\beta,\tilde{X},\alpha).
\]
It attains a finite value $\mathcal{G}^2<+\infty$ if and only if the following conditions are satisfied:
\begin{eqnarray*}
&&1-\hbox{tr}(\tilde{X}N_0)=0,\quad  \beta_l+\hbox{tr}(\tilde{X}M_l)=0, \ l=1,\ldots,d, \quad W-\sum_{j=1}^n\alpha_j\tilde{H}_j=0,
\\
&& A X-b x^T=0, \quad Ax=b, \quad \alpha_j-\hbox{tr}(\tilde{X}N_j)=0, \ j=0,1,\ldots, n.
\end{eqnarray*}
Taking into account the condition $\hbox{tr}(\tilde{X}N_0)=1$ we conclude $\varphi=1$. Hence $\tilde{X}=\left(\begin{array}{cc}
X & x \\
x^T & 1
\end{array}\right)\succeq 0,$ or, equivalently, $X\succeq xx^T$.
The conditions $A X-b x^T=0$ and $Ax=b$ follow from the  identity:
\begin{eqnarray*}
\hbox{tr}(\tilde{X}M_*(u,V))&=&\frac12 \hbox{tr}\left[
\left(\begin{array}{cc}
X & x \\
x^T & 1
\end{array}\right)
\left(\begin{array}{cc}
V^TA+A^TV & A^Tu-V^Tb \\
u^TA-b^TV & -2u^Tb
\end{array}\right)
\right]
\\
&=&\hbox{tr}(XA^TV)-\hbox{tr}(xb^TV)+x^TA^Tu-u^Tb .
\end{eqnarray*}
The conditions $\beta_l+\hbox{tr}(\tilde{X}M_l)=0$ and $\beta\geq 0$ yield $\hbox{tr}(XP_l)+2q_l^Tx+r_l\leq 0.$ for $l=1,\ldots,d$. Since $W\succeq 0$ we deduce $\sum_{j=0}^n\alpha_j\tilde{H}_j\succeq 0$. From the condition $\alpha_j-\hbox{tr}(\tilde{X}N_j)=0$ we obtain $x^Te_j=\alpha_j$, i.~e. $x=(\alpha_1,\ldots, \alpha_n)^T$ and $\alpha_0=1$. Finally, the real LMI $\tilde{H}_0+\sum_{j=1}^n x_j \tilde{H}_j\succeq 0$ is equivalent to the complex LMI ${H}_0+\sum_{j=1}^n x_j {H}_j\succeq 0$, using equivalence \eqref{realvscomplex}. As  $\hbox{tr}(\tilde X M_0) =\hbox{tr}(XP_0)+2q_0^Tx+r_0$ the second Lagrangian dual has the form of SDP \eqref{D3}, as claimed.

\end{document}